\input amstex
%

\def\next{AMS-SEKR}\ifx\styname\next \endinput\fi
\catcode`\@=11
\def\styname{AMS-SEKR}
\def\styversion{2.0}
{\W@{}\W@{\styname.STY - Version \styversion}\W@{}}
\hyphenation{acad-e-my acad-e-mies af-ter-thought anom-aly anom-alies
an-ti-deriv-a-tive an-tin-o-my an-tin-o-mies apoth-e-o-ses apoth-e-o-sis
ap-pen-dix ar-che-typ-al as-sign-a-ble as-sist-ant-ship as-ymp-tot-ic
asyn-chro-nous at-trib-uted at-trib-ut-able bank-rupt bank-rupt-cy
bi-dif-fer-en-tial blue-print busier busiest cat-a-stroph-ic
cat-a-stroph-i-cally con-gress cross-hatched data-base de-fin-i-tive
de-riv-a-tive dis-trib-ute dri-ver dri-vers eco-nom-ics econ-o-mist
elit-ist equi-vari-ant ex-quis-ite ex-tra-or-di-nary flow-chart
for-mi-da-ble forth-right friv-o-lous ge-o-des-ic ge-o-det-ic geo-met-ric
griev-ance griev-ous griev-ous-ly hexa-dec-i-mal ho-lo-no-my ho-mo-thetic
ideals idio-syn-crasy in-fin-ite-ly in-fin-i-tes-i-mal ir-rev-o-ca-ble
key-stroke lam-en-ta-ble light-weight mal-a-prop-ism man-u-script
mar-gin-al meta-bol-ic me-tab-o-lism meta-lan-guage me-trop-o-lis
met-ro-pol-i-tan mi-nut-est mol-e-cule mono-chrome mono-pole mo-nop-oly
mono-spline mo-not-o-nous mul-ti-fac-eted mul-ti-plic-able non-euclid-ean
non-iso-mor-phic non-smooth par-a-digm par-a-bol-ic pa-rab-o-loid
pa-ram-e-trize para-mount pen-ta-gon phe-nom-e-non post-script pre-am-ble
pro-ce-dur-al pro-hib-i-tive pro-hib-i-tive-ly pseu-do-dif-fer-en-tial
pseu-do-fi-nite pseu-do-nym qua-drat-ics quad-ra-ture qua-si-smooth
qua-si-sta-tion-ary qua-si-tri-an-gu-lar quin-tes-sence quin-tes-sen-tial
re-arrange-ment rec-tan-gle ret-ri-bu-tion retro-fit retro-fit-ted
right-eous right-eous-ness ro-bot ro-bot-ics sched-ul-ing se-mes-ter
semi-def-i-nite semi-ho-mo-thet-ic set-up se-vere-ly side-step sov-er-eign
spe-cious spher-oid spher-oid-al star-tling star-tling-ly
sta-tis-tics sto-chas-tic straight-est strange-ness strat-a-gem strong-hold
sum-ma-ble symp-to-matic syn-chro-nous topo-graph-i-cal tra-vers-a-ble
tra-ver-sal tra-ver-sals treach-ery turn-around un-at-tached un-err-ing-ly
white-space wide-spread wing-spread wretch-ed wretch-ed-ly Brown-ian
Eng-lish Euler-ian Feb-ru-ary Gauss-ian Grothen-dieck Hamil-ton-ian
Her-mit-ian Jan-u-ary Japan-ese Kor-te-weg Le-gendre Lip-schitz
Lip-schitz-ian Mar-kov-ian Noe-ther-ian No-vem-ber Rie-mann-ian
Schwarz-schild Sep-tem-ber
form per-iods Uni-ver-si-ty cri-ti-sism for-ma-lism}
\Invalid@\nofrills
\Invalid@\usualspace
\newif\ifnofrills@
\def\nofrills@#1#2{\relaxnext@
  \DN@{\ifx\next\nofrills
    \nofrills@true\let#2\relax\DN@\nofrills{\nextii@}%
  \else
    \nofrills@false\def#2{#1}\let\next@\nextii@\fi
\next@}}
\def\usualspace@#1{\ifnofrills@\def\usualspace{#1}\fi}
\def\addto#1#2{\csname \expandafter\eat@\string#1@\endcsname
  \expandafter{\the\csname \expandafter\eat@\string#1@\endcsname#2}}
\newdimen\bigsize@
\def\big@#1#2{{\hbox{$\left#2\vcenter to#1\bigsize@{}%
  \right.\nulldelimiterspace\z@\m@th$}}}
\def\big{\big@\@ne}
\def\Big{\big@{1.5}}
\def\bigg{\big@\tw@}
\def\Bigg{\big@{2.5}}
\def\raggedcenter@{\leftskip\z@ plus.4\hsize \rightskip\leftskip
 \parfillskip\z@ \parindent\z@ \spaceskip.3333em \xspaceskip.5em
 \pretolerance9999\tolerance9999 \exhyphenpenalty\@M
 \hyphenpenalty\@M \let\\\linebreak}
\def\upperspecialchars{\def\ss{SS}\let\i=I\let\j=J\let\ae\AE\let\oe\OE
  \let\o\O\let\aa\AA\let\l\L}
\def\uppercasetext@#1{%
  {\spaceskip1.2\fontdimen2\the\font plus1.2\fontdimen3\the\font
   \upperspecialchars\uctext@#1$\m@th\aftergroup\eat@$}}
\def\uctext@#1$#2${\endash@#1-\endash@$#2$\uctext@}
\def\endash@#1-#2\endash@{\uppercase{#1}\if\notempty{#2}--\endash@#2\endash@\fi}
\def\runaway@#1{\DN@{#1}\ifx\envir@\next@
  \Err@{You seem to have a missing or misspelled \string\end#1 ...}%
  \let\envir@\empty\fi}
\newif\iftemp@
\def\notempty#1{TT\fi\def\test@{#1}\ifx\test@\empty\temp@false
  \else\temp@true\fi \iftemp@}
\font@\tensmc=cmcsc10
\font@\sevenex=cmex7
\font@\sevenit=cmti7
\font@\eightrm=cmr8 
\font@\sixrm=cmr6 
\font@\eighti=cmmi8     \skewchar\eighti='177 
\font@\sixi=cmmi6       \skewchar\sixi='177   
\font@\eightsy=cmsy8    \skewchar\eightsy='60 
\font@\sixsy=cmsy6      \skewchar\sixsy='60   
\font@\eightex=cmex8
\font@\eightbf=cmbx8 
\font@\sixbf=cmbx6   
\font@\eightit=cmti8 
\font@\eightsl=cmsl8 
\font@\eightsmc=cmcsc8
\font@\eighttt=cmtt8 


\loadmsam
\loadmsbm
\loadeufm
\UseAMSsymbols
\newtoks\tenpoint@
\def\tenpoint{\normalbaselineskip12\p@
 \abovedisplayskip12\p@ plus3\p@ minus9\p@
 \belowdisplayskip\abovedisplayskip
 \abovedisplayshortskip\z@ plus3\p@
 \belowdisplayshortskip7\p@ plus3\p@ minus4\p@
 \textonlyfont@\rm\tenrm \textonlyfont@\it\tenit
 \textonlyfont@\sl\tensl \textonlyfont@\bf\tenbf
 \textonlyfont@\smc\tensmc \textonlyfont@\tt\tentt
 \textonlyfont@\bsmc\tenbsmc
 \ifsyntax@ \def\big##1{{\hbox{$\left##1\right.$}}}%
  \let\Big\big \let\bigg\big \let\Bigg\big
 \else
  \textfont\z@=\tenrm  \scriptfont\z@=\sevenrm  \scriptscriptfont\z@=\fiverm
  \textfont\@ne=\teni  \scriptfont\@ne=\seveni  \scriptscriptfont\@ne=\fivei
  \textfont\tw@=\tensy \scriptfont\tw@=\sevensy \scriptscriptfont\tw@=\fivesy
  \textfont\thr@@=\tenex \scriptfont\thr@@=\sevenex
        \scriptscriptfont\thr@@=\sevenex
  \textfont\itfam=\tenit \scriptfont\itfam=\sevenit
        \scriptscriptfont\itfam=\sevenit
  \textfont\bffam=\tenbf \scriptfont\bffam=\sevenbf
        \scriptscriptfont\bffam=\fivebf
  \setbox\strutbox\hbox{\vrule height8.5\p@ depth3.5\p@ width\z@}%
  \setbox\strutbox@\hbox{\lower.5\normallineskiplimit\vbox{%
        \kern-\normallineskiplimit\copy\strutbox}}%
 \setbox\z@\vbox{\hbox{$($}\kern\z@}\bigsize@=1.2\ht\z@
 \fi
 \normalbaselines\rm\ex@.2326ex\jot3\ex@\the\tenpoint@}
\newtoks\eightpoint@
\def\eightpoint{\normalbaselineskip10\p@
 \abovedisplayskip10\p@ plus2.4\p@ minus7.2\p@
 \belowdisplayskip\abovedisplayskip
 \abovedisplayshortskip\z@ plus2.4\p@
 \belowdisplayshortskip5.6\p@ plus2.4\p@ minus3.2\p@
 \textonlyfont@\rm\eightrm \textonlyfont@\it\eightit
 \textonlyfont@\sl\eightsl \textonlyfont@\bf\eightbf
 \textonlyfont@\smc\eightsmc \textonlyfont@\tt\eighttt
 \textonlyfont@\bsmc\eightbsmc
 \ifsyntax@\def\big##1{{\hbox{$\left##1\right.$}}}%
  \let\Big\big \let\bigg\big \let\Bigg\big
 \else
  \textfont\z@=\eightrm \scriptfont\z@=\sixrm \scriptscriptfont\z@=\fiverm
  \textfont\@ne=\eighti \scriptfont\@ne=\sixi \scriptscriptfont\@ne=\fivei
  \textfont\tw@=\eightsy \scriptfont\tw@=\sixsy \scriptscriptfont\tw@=\fivesy
  \textfont\thr@@=\eightex \scriptfont\thr@@=\sevenex
   \scriptscriptfont\thr@@=\sevenex
  \textfont\itfam=\eightit \scriptfont\itfam=\sevenit
   \scriptscriptfont\itfam=\sevenit
  \textfont\bffam=\eightbf \scriptfont\bffam=\sixbf
   \scriptscriptfont\bffam=\fivebf
 \setbox\strutbox\hbox{\vrule height7\p@ depth3\p@ width\z@}%
 \setbox\strutbox@\hbox{\raise.5\normallineskiplimit\vbox{%
   \kern-\normallineskiplimit\copy\strutbox}}%
 \setbox\z@\vbox{\hbox{$($}\kern\z@}\bigsize@=1.2\ht\z@
 \fi
 \normalbaselines\eightrm\ex@.2326ex\jot3\ex@\the\eightpoint@}

\font@\twelverm=cmr10 scaled\magstep1
\font@\twelveit=cmti10 scaled\magstep1
\font@\twelvesl=cmsl10 scaled\magstep1
\font@\twelvesmc=cmcsc10 scaled\magstep1
\font@\twelvett=cmtt10 scaled\magstep1
\font@\twelvebf=cmbx10 scaled\magstep1
\font@\twelvei=cmmi10 scaled\magstep1
\font@\twelvesy=cmsy10 scaled\magstep1
\font@\twelveex=cmex10 scaled\magstep1
\font@\twelvemsa=msam10 scaled\magstep1
\font@\twelveeufm=eufm10 scaled\magstep1
\font@\twelvemsb=msbm10 scaled\magstep1
\newtoks\twelvepoint@
\def\twelvepoint{\normalbaselineskip15\p@
 \abovedisplayskip15\p@ plus3.6\p@ minus10.8\p@
 \belowdisplayskip\abovedisplayskip
 \abovedisplayshortskip\z@ plus3.6\p@
 \belowdisplayshortskip8.4\p@ plus3.6\p@ minus4.8\p@
 \textonlyfont@\rm\twelverm \textonlyfont@\it\twelveit
 \textonlyfont@\sl\twelvesl \textonlyfont@\bf\twelvebf
 \textonlyfont@\smc\twelvesmc \textonlyfont@\tt\twelvett
 \textonlyfont@\bsmc\twelvebsmc
 \ifsyntax@ \def\big##1{{\hbox{$\left##1\right.$}}}%
  \let\Big\big \let\bigg\big \let\Bigg\big
 \else
  \textfont\z@=\twelverm  \scriptfont\z@=\tenrm  \scriptscriptfont\z@=\sevenrm
  \textfont\@ne=\twelvei  \scriptfont\@ne=\teni  \scriptscriptfont\@ne=\seveni
  \textfont\tw@=\twelvesy \scriptfont\tw@=\tensy \scriptscriptfont\tw@=\sevensy
  \textfont\thr@@=\twelveex \scriptfont\thr@@=\tenex
        \scriptscriptfont\thr@@=\tenex
  \textfont\itfam=\twelveit \scriptfont\itfam=\tenit
        \scriptscriptfont\itfam=\tenit
  \textfont\bffam=\twelvebf \scriptfont\bffam=\tenbf
        \scriptscriptfont\bffam=\sevenbf
  \setbox\strutbox\hbox{\vrule height10.2\p@ depth4.2\p@ width\z@}%
  \setbox\strutbox@\hbox{\lower.6\normallineskiplimit\vbox{%
        \kern-\normallineskiplimit\copy\strutbox}}%
 \setbox\z@\vbox{\hbox{$($}\kern\z@}\bigsize@=1.4\ht\z@
 \fi
 \normalbaselines\rm\ex@.2326ex\jot3.6\ex@\the\twelvepoint@}

\def\headfonts{\twelvepoint\bf}

\font@\fourteenrm=cmr10 scaled\magstep2
\font@\fourteenit=cmti10 scaled\magstep2
\font@\fourteensl=cmsl10 scaled\magstep2
\font@\fourteensmc=cmcsc10 scaled\magstep2
\font@\fourteentt=cmtt10 scaled\magstep2
\font@\fourteenbf=cmbx10 scaled\magstep2
\font@\fourteeni=cmmi10 scaled\magstep2
\font@\fourteensy=cmsy10 scaled\magstep2
\font@\fourteenex=cmex10 scaled\magstep2
\font@\fourteenmsa=msam10 scaled\magstep2
\font@\fourteeneufm=eufm10 scaled\magstep2
\font@\fourteenmsb=msbm10 scaled\magstep2
\newtoks\fourteenpoint@
\def\fourteenpoint{\normalbaselineskip15\p@
 \abovedisplayskip18\p@ plus4.3\p@ minus12.9\p@
 \belowdisplayskip\abovedisplayskip
 \abovedisplayshortskip\z@ plus4.3\p@
 \belowdisplayshortskip10.1\p@ plus4.3\p@ minus5.8\p@
 \textonlyfont@\rm\fourteenrm \textonlyfont@\it\fourteenit
 \textonlyfont@\sl\fourteensl \textonlyfont@\bf\fourteenbf
 \textonlyfont@\smc\fourteensmc \textonlyfont@\tt\fourteentt
 \textonlyfont@\bsmc\fourteenbsmc
 \ifsyntax@ \def\big##1{{\hbox{$\left##1\right.$}}}%
  \let\Big\big \let\bigg\big \let\Bigg\big
 \else
  \textfont\z@=\fourteenrm  \scriptfont\z@=\twelverm  \scriptscriptfont\z@=\tenrm
  \textfont\@ne=\fourteeni  \scriptfont\@ne=\twelvei  \scriptscriptfont\@ne=\teni
  \textfont\tw@=\fourteensy \scriptfont\tw@=\twelvesy \scriptscriptfont\tw@=\tensy
  \textfont\thr@@=\fourteenex \scriptfont\thr@@=\twelveex
        \scriptscriptfont\thr@@=\twelveex
  \textfont\itfam=\fourteenit \scriptfont\itfam=\twelveit
        \scriptscriptfont\itfam=\twelveit
  \textfont\bffam=\fourteenbf \scriptfont\bffam=\twelvebf
        \scriptscriptfont\bffam=\tenbf
  \setbox\strutbox\hbox{\vrule height12.2\p@ depth5\p@ width\z@}%
  \setbox\strutbox@\hbox{\lower.72\normallineskiplimit\vbox{%
        \kern-\normallineskiplimit\copy\strutbox}}%
 \setbox\z@\vbox{\hbox{$($}\kern\z@}\bigsize@=1.7\ht\z@
 \fi
 \normalbaselines\rm\ex@.2326ex\jot4.3\ex@\the\fourteenpoint@}

\def\chapheadfonts{\fourteenpoint\bf}

\font@\seventeenrm=cmr10 scaled\magstep3
\font@\seventeenit=cmti10 scaled\magstep3
\font@\seventeensl=cmsl10 scaled\magstep3
\font@\seventeensmc=cmcsc10 scaled\magstep3
\font@\seventeentt=cmtt10 scaled\magstep3
\font@\seventeenbf=cmbx10 scaled\magstep3
\font@\seventeeni=cmmi10 scaled\magstep3
\font@\seventeensy=cmsy10 scaled\magstep3
\font@\seventeenex=cmex10 scaled\magstep3
\font@\seventeenmsa=msam10 scaled\magstep3
\font@\seventeeneufm=eufm10 scaled\magstep3
\font@\seventeenmsb=msbm10 scaled\magstep3
\newtoks\seventeenpoint@
\def\seventeenpoint{\normalbaselineskip18\p@
 \abovedisplayskip21.6\p@ plus5.2\p@ minus15.4\p@
 \belowdisplayskip\abovedisplayskip
 \abovedisplayshortskip\z@ plus5.2\p@
 \belowdisplayshortskip12.1\p@ plus5.2\p@ minus7\p@
 \textonlyfont@\rm\seventeenrm \textonlyfont@\it\seventeenit
 \textonlyfont@\sl\seventeensl \textonlyfont@\bf\seventeenbf
 \textonlyfont@\smc\seventeensmc \textonlyfont@\tt\seventeentt
 \textonlyfont@\bsmc\seventeenbsmc
 \ifsyntax@ \def\big##1{{\hbox{$\left##1\right.$}}}%
  \let\Big\big \let\bigg\big \let\Bigg\big
 \else
  \textfont\z@=\seventeenrm  \scriptfont\z@=\fourteenrm  \scriptscriptfont\z@=\twelverm
  \textfont\@ne=\seventeeni  \scriptfont\@ne=\fourteeni  \scriptscriptfont\@ne=\twelvei
  \textfont\tw@=\seventeensy \scriptfont\tw@=\fourteensy \scriptscriptfont\tw@=\twelvesy
  \textfont\thr@@=\seventeenex \scriptfont\thr@@=\fourteenex
        \scriptscriptfont\thr@@=\fourteenex
  \textfont\itfam=\seventeenit \scriptfont\itfam=\fourteenit
        \scriptscriptfont\itfam=\fourteenit
  \textfont\bffam=\seventeenbf \scriptfont\bffam=\fourteenbf
        \scriptscriptfont\bffam=\twelvebf
  \setbox\strutbox\hbox{\vrule height14.6\p@ depth6\p@ width\z@}%
  \setbox\strutbox@\hbox{\lower.86\normallineskiplimit\vbox{%
        \kern-\normallineskiplimit\copy\strutbox}}%
 \setbox\z@\vbox{\hbox{$($}\kern\z@}\bigsize@=2\ht\z@
 \fi
 \normalbaselines\rm\ex@.2326ex\jot5.2\ex@\the\seventeenpoint@}

\font@\rrrrrm=cmr10 scaled\magstep4
\font@\bigtitlefont=cmbx10 scaled\magstep4

\parindent1pc
\normallineskiplimit\p@
\newdimen\indenti \indenti=2pc
\def\pageheight#1{\vsize#1}
\def\pagewidth#1{\hsize#1%
   \captionwidth@\hsize \advance\captionwidth@-2\indenti}
\pagewidth{30pc} \pageheight{47pc}
\def\topmatter{%
 \ifx\undefined\msafam
 \else\font@\eightmsa=msam8 \font@\sixmsa=msam6
   \ifsyntax@\else \addto\tenpoint{\textfont\msafam=\tenmsa
              \scriptfont\msafam=\sevenmsa \scriptscriptfont\msafam=\fivemsa}%
     \addto\eightpoint{\textfont\msafam=\eightmsa \scriptfont\msafam=\sixmsa
              \scriptscriptfont\msafam=\fivemsa}%
   \fi
 \fi
 \ifx\undefined\msbfam
 \else\font@\eightmsb=msbm8 \font@\sixmsb=msbm6
   \ifsyntax@\else \addto\tenpoint{\textfont\msbfam=\tenmsb
         \scriptfont\msbfam=\sevenmsb \scriptscriptfont\msbfam=\fivemsb}%
     \addto\eightpoint{\textfont\msbfam=\eightmsb \scriptfont\msbfam=\sixmsb
         \scriptscriptfont\msbfam=\fivemsb}%
   \fi
 \fi
 \ifx\undefined\eufmfam
 \else \font@\eighteufm=eufm8 \font@\sixeufm=eufm6
   \ifsyntax@\else \addto\tenpoint{\textfont\eufmfam=\teneufm
       \scriptfont\eufmfam=\seveneufm \scriptscriptfont\eufmfam=\fiveeufm}%
     \addto\eightpoint{\textfont\eufmfam=\eighteufm
       \scriptfont\eufmfam=\sixeufm \scriptscriptfont\eufmfam=\fiveeufm}%
   \fi
 \fi
 \ifx\undefined\eufbfam
 \else \font@\eighteufb=eufb8 \font@\sixeufb=eufb6
   \ifsyntax@\else \addto\tenpoint{\textfont\eufbfam=\teneufb
      \scriptfont\eufbfam=\seveneufb \scriptscriptfont\eufbfam=\fiveeufb}%
    \addto\eightpoint{\textfont\eufbfam=\eighteufb
      \scriptfont\eufbfam=\sixeufb \scriptscriptfont\eufbfam=\fiveeufb}%
   \fi
 \fi
 \ifx\undefined\eusmfam
 \else \font@\eighteusm=eusm8 \font@\sixeusm=eusm6
   \ifsyntax@\else \addto\tenpoint{\textfont\eusmfam=\teneusm
       \scriptfont\eusmfam=\seveneusm \scriptscriptfont\eusmfam=\fiveeusm}%
     \addto\eightpoint{\textfont\eusmfam=\eighteusm
       \scriptfont\eusmfam=\sixeusm \scriptscriptfont\eusmfam=\fiveeusm}%
   \fi
 \fi
 \ifx\undefined\eusbfam
 \else \font@\eighteusb=eusb8 \font@\sixeusb=eusb6
   \ifsyntax@\else \addto\tenpoint{\textfont\eusbfam=\teneusb
       \scriptfont\eusbfam=\seveneusb \scriptscriptfont\eusbfam=\fiveeusb}%
     \addto\eightpoint{\textfont\eusbfam=\eighteusb
       \scriptfont\eusbfam=\sixeusb \scriptscriptfont\eusbfam=\fiveeusb}%
   \fi
 \fi
 \ifx\undefined\eurmfam
 \else \font@\eighteurm=eurm8 \font@\sixeurm=eurm6
   \ifsyntax@\else \addto\tenpoint{\textfont\eurmfam=\teneurm
       \scriptfont\eurmfam=\seveneurm \scriptscriptfont\eurmfam=\fiveeurm}%
     \addto\eightpoint{\textfont\eurmfam=\eighteurm
       \scriptfont\eurmfam=\sixeurm \scriptscriptfont\eurmfam=\fiveeurm}%
   \fi
 \fi
 \ifx\undefined\eurbfam
 \else \font@\eighteurb=eurb8 \font@\sixeurb=eurb6
   \ifsyntax@\else \addto\tenpoint{\textfont\eurbfam=\teneurb
       \scriptfont\eurbfam=\seveneurb \scriptscriptfont\eurbfam=\fiveeurb}%
    \addto\eightpoint{\textfont\eurbfam=\eighteurb
       \scriptfont\eurbfam=\sixeurb \scriptscriptfont\eurbfam=\fiveeurb}%
   \fi
 \fi
 \ifx\undefined\cmmibfam
 \else \font@\eightcmmib=cmmib8 \font@\sixcmmib=cmmib6
   \ifsyntax@\else \addto\tenpoint{\textfont\cmmibfam=\tencmmib
       \scriptfont\cmmibfam=\sevencmmib \scriptscriptfont\cmmibfam=\fivecmmib}%
    \addto\eightpoint{\textfont\cmmibfam=\eightcmmib
       \scriptfont\cmmibfam=\sixcmmib \scriptscriptfont\cmmibfam=\fivecmmib}%
   \fi
 \fi
 \ifx\undefined\cmbsyfam
 \else \font@\eightcmbsy=cmbsy8 \font@\sixcmbsy=cmbsy6
   \ifsyntax@\else \addto\tenpoint{\textfont\cmbsyfam=\tencmbsy
      \scriptfont\cmbsyfam=\sevencmbsy \scriptscriptfont\cmbsyfam=\fivecmbsy}%
    \addto\eightpoint{\textfont\cmbsyfam=\eightcmbsy
      \scriptfont\cmbsyfam=\sixcmbsy \scriptscriptfont\cmbsyfam=\fivecmbsy}%
   \fi
 \fi
 \let\topmatter\relax}
\def\chapterno@{\uppercase\expandafter{\romannumeral\chaptercount@}}
\newcount\chaptercount@
\def\chapter{\nofrills@{\afterassignment\chapterno@
                        CHAPTER \global\chaptercount@=}\chapter@
 \DNii@##1{\leavevmode\hskip-\leftskip
   \rlap{\vbox to\z@{\vss\centerline{\eightpoint
   \chapter@##1\unskip}\baselineskip2pc\null}}\hskip\leftskip
   \nofrills@false}%
 \FN@\next@}
\newbox\titlebox@

\def\title{\nofrills@{\relax}\title@%
 \DNii@##1\endtitle{\global\setbox\titlebox@\vtop{\tenpoint\bf
 \raggedcenter@\ignorespaces
 \baselineskip1.3\baselineskip\title@{##1}\endgraf}%
 \ifmonograph@ \edef\next{\the\leftheadtoks}\ifx\next\empty
    \leftheadtext{##1}\fi
 \fi
 \edef\next{\the\rightheadtoks}\ifx\next\empty \rightheadtext{##1}\fi
 }\FN@\next@}
\newbox\authorbox@
\def\author#1\endauthor{\global\setbox\authorbox@
 \vbox{\tenpoint\smc\raggedcenter@\ignorespaces
 #1\endgraf}\relaxnext@ \edef\next{\the\leftheadtoks}%
 \ifx\next\empty\leftheadtext{#1}\fi}
\newbox\affilbox@
\def\affil#1\endaffil{\global\setbox\affilbox@
 \vbox{\tenpoint\raggedcenter@\ignorespaces#1\endgraf}}
\newcount\addresscount@
\addresscount@\z@
\def\address#1\endaddress{\global\advance\addresscount@\@ne
  \expandafter\gdef\csname address\number\addresscount@\endcsname
  {\vskip12\p@ minus6\p@\noindent\eightpoint\smc\ignorespaces#1\par}}
\def\email{\nofrills@{\eightpoint{\it E-mail\/}:\enspace}\email@
  \DNii@##1\endemail{%
  \expandafter\gdef\csname email\number\addresscount@\endcsname
  {\def\usualspace{{\it\enspace}}\smallskip\noindent\eightpoint\email@
  \ignorespaces##1\par}}%
 \FN@\next@}
\def\thedate@{}
\def\date#1\enddate{\gdef\thedate@{\tenpoint\ignorespaces#1\unskip}}
\def\thethanks@{}
\def\thanks#1\endthanks{\gdef\thethanks@{\eightpoint\ignorespaces#1.\unskip}}
\def\thekeywords@{}
\def\keywords{\nofrills@{{\it Key words and phrases.\enspace}}\keywords@
 \DNii@##1\endkeywords{\def\thekeywords@{\def\usualspace{{\it\enspace}}%
 \eightpoint\keywords@\ignorespaces##1\unskip.}}%
 \FN@\next@}
\def\thesubjclass@{}
\def\subjclass{\nofrills@{{\rm2010 {\it Mathematics Subject
   Classification\/}.\enspace}}\subjclass@
 \DNii@##1\endsubjclass{\def\thesubjclass@{\def\usualspace
  {{\rm\enspace}}\eightpoint\subjclass@\ignorespaces##1\unskip.}}%
 \FN@\next@}
\newbox\abstractbox@
\def\abstract{\nofrills@{{\smc Abstract.\enspace}}\abstract@
 \DNii@{\setbox\abstractbox@\vbox\bgroup\noindent$$\vbox\bgroup
  \def\envir@{abstract}\advance\hsize-2\indenti
  \usualspace@{{\enspace}}\eightpoint \noindent\abstract@\ignorespaces}%
 \FN@\next@}
\def\endabstract{\par\unskip\egroup$$\egroup}
\def\widestnumber#1#2{\begingroup\let\head\null\let\subhead\empty
   \let\subsubhead\subhead
   \ifx#1\head\global\setbox\tocheadbox@\hbox{#2.\enspace}%
   \else\ifx#1\subhead\global\setbox\tocsubheadbox@\hbox{#2.\enspace}%
   \else\ifx#1\key\bgroup\let\endrefitem@\egroup
        \key#2\endrefitem@\global\refindentwd\wd\keybox@
   \else\ifx#1\no\bgroup\let\endrefitem@\egroup
        \no#2\endrefitem@\global\refindentwd\wd\nobox@
   \else\ifx#1\page\global\setbox\pagesbox@\hbox{\quad\bf#2}%
   \else\ifx#1\item\setboxz@h{#2}\global\rosteritemwd\wdz@
        \global\advance\rosteritemwd by.5\parindent
   \else\message{\string\widestnumber is not defined for this option
   (\string#1)}%
\fi\fi\fi\fi\fi\fi\endgroup}
\newif\ifmonograph@
\def\Monograph{\monograph@true \let\headmark\rightheadtext
  \let\varindent@\indent \def\headfont@{\bf}\def\proclaimheadfont@{\smc}%
  \def\demofont@{\smc}}
\let\varindent@\indent

\newbox\tocheadbox@    \newbox\tocsubheadbox@
\newbox\tocbox@
\def\toc{\toc@{Contents}}
\def\newtocdefs{%
   \def \title##1\endtitle
       {\penaltyandskip@\z@\smallskipamount
        \hangindent\wd\tocheadbox@\noindent{\bf##1}}%
   \def \chapter##1{%
        Chapter \uppercase\expandafter{\romannumeral##1.\unskip}\enspace}%
   \def \specialhead##1\endspecialhead
       {\par\hangindent\wd\tocheadbox@ \noindent##1\par}%
   \def \head##1 ##2\endhead
       {\par\hangindent\wd\tocheadbox@ \noindent
        \if\notempty{##1}\hbox to\wd\tocheadbox@{\hfil##1\enspace}\fi
        ##2\par}%
   \def \subhead##1 ##2\endsubhead
       {\par\vskip-\parskip {\normalbaselines
        \advance\leftskip\wd\tocheadbox@
        \hangindent\wd\tocsubheadbox@ \noindent
        \if\notempty{##1}\hbox to\wd\tocsubheadbox@{##1\unskip\hfil}\fi
         ##2\par}}%
   \def \subsubhead##1 ##2\endsubsubhead
       {\par\vskip-\parskip {\normalbaselines
        \advance\leftskip\wd\tocheadbox@
        \hangindent\wd\tocsubheadbox@ \noindent
        \if\notempty{##1}\hbox to\wd\tocsubheadbox@{##1\unskip\hfil}\fi
        ##2\par}}}
\def\toc@#1{\relaxnext@
   \def\page##1%
       {\unskip\penalty0\null\hfil
        \rlap{\hbox to\wd\pagesbox@{\quad\hfil##1}}\hfilneg\penalty\@M}%
 \DN@{\ifx\next\nofrills\DN@\nofrills{\nextii@}%
      \else\DN@{\nextii@{{#1}}}\fi
      \next@}%
 \DNii@##1{%
\ifmonograph@\bgroup\else\setbox\tocbox@\vbox\bgroup
   \centerline{\headfont@\ignorespaces##1\unskip}\nobreak
   \vskip\belowheadskip \fi
   \setbox\tocheadbox@\hbox{0.\enspace}%
   \setbox\tocsubheadbox@\hbox{0.0.\enspace}%
   \leftskip\indenti \rightskip\leftskip
   \setbox\pagesbox@\hbox{\bf\quad000}\advance\rightskip\wd\pagesbox@
   \newtocdefs
 }%
 \FN@\next@}
\def\endtoc{\par\egroup}
\let\pretitle\relax
\let\preauthor\relax
\let\preaffil\relax
\let\predate\relax
\let\preabstract\relax
\let\prepaper\relax
\def\dedicatory #1\enddedicatory{\def\preabstract{{\medskip
  \eightpoint\it \raggedcenter@#1\endgraf}}}
\def\thetranslator@{}
\def\translator#1\endtranslator{\def\thetranslator@{\nobreak\medskip
 \line{\eightpoint\hfil Translated by \uppercase{#1}\qquad\qquad}\nobreak}}
\outer\def\endtopmatter{\runaway@{abstract}%
 \edef\next{\the\leftheadtoks}\ifx\next\empty
  \expandafter\leftheadtext\expandafter{\the\rightheadtoks}\fi
 \ifmonograph@\else
   \ifx\thesubjclass@\empty\else \makefootnote@{}{\thesubjclass@}\fi
   \ifx\thekeywords@\empty\else \makefootnote@{}{\thekeywords@}\fi
   \ifx\thethanks@\empty\else \makefootnote@{}{\thethanks@}\fi
 \fi
  \pretitle
  \ifmonograph@ \topskip7pc \else \topskip4pc \fi
  \box\titlebox@
  \topskip10pt
  \preauthor
  \ifvoid\authorbox@\else \vskip2.5pc plus1pc \unvbox\authorbox@\fi
  \preaffil
  \ifvoid\affilbox@\else \vskip1pc plus.5pc \unvbox\affilbox@\fi
  \predate
  \ifx\thedate@\empty\else \vskip1pc plus.5pc \line{\hfil\thedate@\hfil}\fi
  \preabstract
  \ifvoid\abstractbox@\else \vskip1.5pc plus.5pc \unvbox\abstractbox@ \fi
  \ifvoid\tocbox@\else\vskip1.5pc plus.5pc \unvbox\tocbox@\fi
  \prepaper
  \vskip2pc plus1pc
}
\def\document{\let\fontlist@\relax\let\alloclist@\relax
  \tenpoint}

\newskip\aboveheadskip       \aboveheadskip1.8\bigskipamount
\newdimen\belowheadskip      \belowheadskip1.8\medskipamount

\def\headfont@{\smc}
\def\penaltyandskip@#1#2{\relax\ifdim\lastskip<#2\relax\removelastskip
      \ifnum#1=\z@\else\penalty@#1\relax\fi\vskip#2%
  \else\ifnum#1=\z@\else\penalty@#1\relax\fi\fi}
\def\nobreak{\penalty\@M
  \ifvmode\def\penalty@{\let\penalty@\penalty\count@@@}%
  \everypar{\let\penalty@\penalty\everypar{}}\fi}
\let\penalty@\penalty
\def\heading#1\endheading{\head#1\endhead}

\def\specialheadfont@{\bf}
\outer\def\specialhead{\par\penaltyandskip@{-200}\aboveheadskip
  \begingroup\interlinepenalty\@M\rightskip\z@ plus\hsize \let\\\linebreak
  \specialheadfont@\noindent\ignorespaces}
\def\endspecialhead{\par\endgroup\nobreak\vskip\belowheadskip}
\let\headmark\eat@
\newskip\subheadskip       \subheadskip\medskipamount
\def\subheadfont@{\bf}
\outer\def\subhead{\nofrills@{.\enspace}\subhead@
 \DNii@##1\endsubhead{\par\penaltyandskip@{-100}\subheadskip
  \varindent@{\usualspace@{{\subheadfont@\enspace}}%
 \subheadfont@\ignorespaces##1\unskip\subhead@}\ignorespaces}%
 \FN@\next@}
\outer\def\subsubhead{\nofrills@{.\enspace}\subsubhead@
 \DNii@##1\endsubsubhead{\par\penaltyandskip@{-50}\medskipamount
      {\usualspace@{{\it\enspace}}%
  \it\ignorespaces##1\unskip\subsubhead@}\ignorespaces}%
 \FN@\next@}
\def\proclaimheadfont@{\bf}
\outer\def\proclaim{\runaway@{proclaim}\def\envir@{proclaim}%
  \nofrills@{.\enspace}\proclaim@
 \DNii@##1{\penaltyandskip@{-100}\medskipamount\varindent@
   \usualspace@{{\proclaimheadfont@\enspace}}\proclaimheadfont@
   \ignorespaces##1\unskip\proclaim@
  \sl\ignorespaces}%
 \FN@\next@}
\outer\def\endproclaim{\let\envir@\relax\par\rm
  \penaltyandskip@{55}\medskipamount}
\def\demoheadfont@{\it}
\def\demo{\runaway@{proclaim}\nofrills@{.\enspace}\demo@
     \DNii@##1{\par\penaltyandskip@\z@\medskipamount
  {\usualspace@{{\demoheadfont@\enspace}}%
  \varindent@\demoheadfont@\ignorespaces##1\unskip\demo@}\rm
  \ignorespaces}\FN@\next@}
\def\enddemo{\par\medskip}
\def\qed{\ifhmode\unskip\nobreak\fi\quad\ifmmode\square\else$\m@th\square$\fi}
\let\remark\demo
\let\endremark\enddemo
\def\definition{\runaway@{proclaim}%
  \nofrills@{.\demoheadfont@\enspace}\definition@
        \DNii@##1{\penaltyandskip@{-100}\medskipamount
        {\usualspace@{{\demoheadfont@\enspace}}%
        \varindent@\demoheadfont@\ignorespaces##1\unskip\definition@}%
        \rm \ignorespaces}\FN@\next@}


\newdimen\rosteritemwd
\newcount\rostercount@
\newif\iffirstitem@
\let\plainitem@\item
\newtoks\everypartoks@
\def\par@{\everypartoks@\expandafter{\the\everypar}\everypar{}}
\def\roster{\edef\leftskip@{\leftskip\the\leftskip}%
 \relaxnext@
 \rostercount@\z@  
 \def\item{\FN@\rosteritem@}%
 \DN@{\ifx\next\runinitem\let\next@\nextii@\else
  \let\next@\nextiii@\fi\next@}%
 \DNii@\runinitem  
  {\unskip  
   \DN@{\ifx\next[\let\next@\nextii@\else
    \ifx\next"\let\next@\nextiii@\else\let\next@\nextiv@\fi\fi\next@}%
   \DNii@[####1]{\rostercount@####1\relax
    \enspace{\rm(\number\rostercount@)}~\ignorespaces}%
   \def\nextiii@"####1"{\enspace{\rm####1}~\ignorespaces}%
   \def\nextiv@{\enspace{\rm(1)}\rostercount@\@ne~}%
   \par@\firstitem@false  
   \FN@\next@}%
 \def\nextiii@{\par\par@  
  \penalty\@m\smallskip\vskip-\parskip
  \firstitem@true}%
 \FN@\next@}
\def\rosteritem@{\iffirstitem@\firstitem@false\else\par\vskip-\parskip\fi
 \leftskip3\parindent\noindent  
 \DNii@[##1]{\rostercount@##1\relax
  \llap{\hbox to2.5\parindent{\hss\rm(\number\rostercount@)}%
   \hskip.5\parindent}\ignorespaces}%
 \def\nextiii@"##1"{%
  \llap{\hbox to2.5\parindent{\hss\rm##1}\hskip.5\parindent}\ignorespaces}%
 \def\nextiv@{\advance\rostercount@\@ne
  \llap{\hbox to2.5\parindent{\hss\rm(\number\rostercount@)}%
   \hskip.5\parindent}}%
 \ifx\next[\let\next@\nextii@\else\ifx\next"\let\next@\nextiii@\else
  \let\next@\nextiv@\fi\fi\next@}

\newif\ifnextRunin@
\def\endroster{\relaxnext@
 \par\leftskip@  
 \penalty-50 \vskip-\parskip\smallskip  
 \DN@{\ifx\next\Runinitem\let\next@\relax
  \else\nextRunin@false\let\item\plainitem@  
   \ifx\next\par 
    \DN@\par{\everypar\expandafter{\the\everypartoks@}}%
   \else  
    \DN@{\noindent\everypar\expandafter{\the\everypartoks@}}%
  \fi\fi\next@}%
 \FN@\next@}
\newcount\rosterhangafter@
\def\Runinitem#1\roster\runinitem{\relaxnext@
 \rostercount@\z@ 
 \def\item{\FN@\rosteritem@}%
 \def\runinitem@{#1}%
 \DN@{\ifx\next[\let\next\nextii@\else\ifx\next"\let\next\nextiii@
  \else\let\next\nextiv@\fi\fi\next}%
 \DNii@[##1]{\rostercount@##1\relax
  \def\item@{{\rm(\number\rostercount@)}}\nextv@}%
 \def\nextiii@"##1"{\def\item@{{\rm##1}}\nextv@}%
 \def\nextiv@{\advance\rostercount@\@ne
  \def\item@{{\rm(\number\rostercount@)}}\nextv@}%
 \def\nextv@{\setbox\z@\vbox  
  {\ifnextRunin@\noindent\fi  
  \runinitem@\unskip\enspace\item@~\par  
  \global\rosterhangafter@\prevgraf}%
  \firstitem@false  
  \ifnextRunin@\else\par\fi  
  \hangafter\rosterhangafter@\hangindent3\parindent
  \ifnextRunin@\noindent\fi  
  \runinitem@\unskip\enspace 
  \item@~\ifnextRunin@\else\par@\fi  
  \nextRunin@true\ignorespaces}%
 \FN@\next@}
\def\footmarkform@#1{$\m@th^{#1}$}
\let\thefootnotemark\footmarkform@
\def\makefootnote@#1#2{\insert\footins
 {\interlinepenalty\interfootnotelinepenalty
 \eightpoint\splittopskip\ht\strutbox\splitmaxdepth\dp\strutbox
 \floatingpenalty\@MM\leftskip\z@\rightskip\z@\spaceskip\z@\xspaceskip\z@
 \leavevmode{#1}\footstrut\ignorespaces#2\unskip\lower\dp\strutbox
 \vbox to\dp\strutbox{}}}
\newcount\footmarkcount@
\footmarkcount@\z@
\def\footnotemark{\let\@sf\empty\relaxnext@
 \ifhmode\edef\@sf{\spacefactor\the\spacefactor}\/\fi
 \DN@{\ifx[\next\let\next@\nextii@\else
  \ifx"\next\let\next@\nextiii@\else
  \let\next@\nextiv@\fi\fi\next@}%
 \DNii@[##1]{\footmarkform@{##1}\@sf}%
 \def\nextiii@"##1"{{##1}\@sf}%
 \def\nextiv@{\iffirstchoice@\global\advance\footmarkcount@\@ne\fi
  \footmarkform@{\number\footmarkcount@}\@sf}%
 \FN@\next@}
\def\footnotetext{\relaxnext@
 \DN@{\ifx[\next\let\next@\nextii@\else
  \ifx"\next\let\next@\nextiii@\else
  \let\next@\nextiv@\fi\fi\next@}%
 \DNii@[##1]##2{\makefootnote@{\footmarkform@{##1}}{##2}}%
 \def\nextiii@"##1"##2{\makefootnote@{##1}{##2}}%
 \def\nextiv@##1{\makefootnote@{\footmarkform@{\number\footmarkcount@}}{##1}}%
 \FN@\next@}
\def\footnote{\let\@sf\empty\relaxnext@
 \ifhmode\edef\@sf{\spacefactor\the\spacefactor}\/\fi
 \DN@{\ifx[\next\let\next@\nextii@\else
  \ifx"\next\let\next@\nextiii@\else
  \let\next@\nextiv@\fi\fi\next@}%
 \DNii@[##1]##2{\footnotemark[##1]\footnotetext[##1]{##2}}%
 \def\nextiii@"##1"##2{\footnotemark"##1"\footnotetext"##1"{##2}}%
 \def\nextiv@##1{\footnotemark\footnotetext{##1}}%
 \FN@\next@}
\def\adjustfootnotemark#1{\advance\footmarkcount@#1\relax}
\def\footnoterule{\kern-3\p@
  \hrule width 5pc\kern 2.6\p@} 
\def\captionfont@{\smc}
\def\topcaption#1#2\endcaption{%
  {\dimen@\hsize \advance\dimen@-\captionwidth@
   \rm\raggedcenter@ \advance\leftskip.5\dimen@ \rightskip\leftskip
  {\captionfont@#1}%
  \if\notempty{#2}.\enspace\ignorespaces#2\fi
  \endgraf}\nobreak\bigskip}
\def\botcaption#1#2\endcaption{%
  \nobreak\bigskip
  \setboxz@h{\captionfont@#1\if\notempty{#2}.\enspace\rm#2\fi}%
  {\dimen@\hsize \advance\dimen@-\captionwidth@
   \leftskip.5\dimen@ \rightskip\leftskip
   \noindent \ifdim\wdz@>\captionwidth@ 
   \else\hfil\fi 
  {\captionfont@#1}\if\notempty{#2}.\enspace\rm#2\fi\endgraf}}
\def\@ins{\par\begingroup\def\vspace##1{\vskip##1\relax}%
  \def\captionwidth##1{\captionwidth@##1\relax}%
  \setbox\z@\vbox\bgroup} 
\def\block{\RIfMIfI@\nondmatherr@\block\fi
       \else\ifvmode\vskip\abovedisplayskip\noindent\fi
        $$\def\endblock{\par\egroup$$}\fi
  \vbox\bgroup\advance\hsize-2\indenti\noindent}
\def\endblock{\par\egroup}
\def\cite#1{{\rm[{\citefont@\m@th#1}]}}
\def\citefont@{\rm}
\def\refsfont@{\eightpoint}
\outer\def\Refs{\runaway@{proclaim}%
 \relaxnext@ \DN@{\ifx\next\nofrills\DN@\nofrills{\nextii@}\else
  \DN@{\nextii@{References}}\fi\next@}%
 \DNii@##1{\penaltyandskip@{-200}\aboveheadskip
  \line{\hfil\headfont@\ignorespaces##1\unskip\hfil}\nobreak
  \vskip\belowheadskip
  \begingroup\refsfont@\sfcode`.=\@m}%
 \FN@\next@}
\def\endRefs{\par\endgroup}
\newbox\nobox@            \newbox\keybox@           \newbox\bybox@
\newbox\paperbox@         \newbox\paperinfobox@     \newbox\jourbox@
\newbox\volbox@           \newbox\issuebox@         \newbox\yrbox@
\newbox\pagesbox@         \newbox\bookbox@          \newbox\bookinfobox@
\newbox\publbox@          \newbox\publaddrbox@      \newbox\finalinfobox@
\newbox\edsbox@           \newbox\langbox@
\newif\iffirstref@        \newif\iflastref@
\newif\ifprevjour@        \newif\ifbook@            \newif\ifprevinbook@
\newif\ifquotes@          \newif\ifbookquotes@      \newif\ifpaperquotes@
\newdimen\bysamerulewd@
\setboxz@h{\refsfont@\kern3em}
\bysamerulewd@\wdz@
\newdimen\refindentwd
\setboxz@h{\refsfont@ 00. }
\refindentwd\wdz@
\outer\def\ref{\begingroup \noindent\hangindent\refindentwd
 \firstref@true \def\nofrills{\def\refkern@{\kern3sp}}%
 \ref@}
\def\ref@{\book@false \bgroup\let\endrefitem@\egroup \ignorespaces}
\def\moreref{\endrefitem@\endref@\firstref@false\ref@}%
\def\transl{\endrefitem@\endref@\firstref@false
  \book@false
  \prepunct@
  \setboxz@h\bgroup \aftergroup\unhbox\aftergroup\z@
    \def\endrefitem@{\unskip\refkern@\egroup}\ignorespaces}%
\def\emptyifempty@{\dimen@\wd\currbox@
  \advance\dimen@-\wd\z@ \advance\dimen@-.1\p@
  \ifdim\dimen@<\z@ \setbox\currbox@\copy\voidb@x \fi}
\let\refkern@\relax
\def\endrefitem@{\unskip\refkern@\egroup
  \setboxz@h{\refkern@}\emptyifempty@}\ignorespaces
\def\refdef@#1#2#3{\edef\next@{\noexpand\endrefitem@
  \let\noexpand\currbox@\csname\expandafter\eat@\string#1box@\endcsname
    \noexpand\setbox\noexpand\currbox@\hbox\bgroup}%
  \toks@\expandafter{\next@}%
  \if\notempty{#2#3}\toks@\expandafter{\the\toks@
  \def\endrefitem@{\unskip#3\refkern@\egroup
  \setboxz@h{#2#3\refkern@}\emptyifempty@}#2}\fi
  \toks@\expandafter{\the\toks@\ignorespaces}%
  \edef#1{\the\toks@}}
\refdef@\no{}{. }
\refdef@\key{[\m@th}{] }
\refdef@\by{}{}
\def\bysame{\by\hbox to\bysamerulewd@{\hrulefill}\thinspace
   \kern0sp}
\def\manyby{\message{\string\manyby is no longer necessary; \string\by
  can be used instead, starting with version 2.0 of \styname.STY}\by}
\refdef@\paper{\ifpaperquotes@``\fi\it}{}
\refdef@\paperinfo{}{}
\def\jour{\endrefitem@\let\currbox@\jourbox@
  \setbox\currbox@\hbox\bgroup
  \def\endrefitem@{\unskip\refkern@\egroup
    \setboxz@h{\refkern@}\emptyifempty@
    \ifvoid\jourbox@\else\prevjour@true\fi}%
\ignorespaces}
\refdef@\vol{\ifbook@\else\bf\fi}{}
\refdef@\issue{no. }{}
\refdef@\yr{}{}
\refdef@\pages{}{}
\def\page{\endrefitem@\def\pp@{\def\pp@{pp.~}p.~}\let\currbox@\pagesbox@
  \setbox\currbox@\hbox\bgroup\ignorespaces}
\def\pp@{pp.~}
\def\book{\endrefitem@ \let\currbox@\bookbox@
 \setbox\currbox@\hbox\bgroup\def\endrefitem@{\unskip\refkern@\egroup
  \setboxz@h{\ifbookquotes@``\fi}\emptyifempty@
  \ifvoid\bookbox@\else\book@true\fi}%
  \ifbookquotes@``\fi\it\ignorespaces}
\def\inbook{\endrefitem@
  \let\currbox@\bookbox@\setbox\currbox@\hbox\bgroup
  \def\endrefitem@{\unskip\refkern@\egroup
  \setboxz@h{\ifbookquotes@``\fi}\emptyifempty@
  \ifvoid\bookbox@\else\book@true\previnbook@true\fi}%
  \ifbookquotes@``\fi\ignorespaces}
\refdef@\eds{(}{, eds.)}
\def\ed{\endrefitem@\let\currbox@\edsbox@
 \setbox\currbox@\hbox\bgroup
 \def\endrefitem@{\unskip, ed.)\refkern@\egroup
  \setboxz@h{(, ed.)}\emptyifempty@}(\ignorespaces}
\refdef@\bookinfo{}{}
\refdef@\publ{}{}
\refdef@\publaddr{}{}
\refdef@\finalinfo{}{}
\refdef@\lang{(}{)}

\let\refdef@\relax 
\def\ppunbox@#1{\ifvoid#1\else\prepunct@\unhbox#1\fi}
\def\nocomma@#1{\ifvoid#1\else\changepunct@3\prepunct@\unhbox#1\fi}
\def\changepunct@#1{\ifnum\lastkern<3 \unkern\kern#1sp\fi}
\def\prepunct@{\count@\lastkern\unkern
  \ifnum\lastpenalty=0
    \let\penalty@\relax
  \else
    \edef\penalty@{\penalty\the\lastpenalty\relax}%
  \fi
  \unpenalty
  \let\refspace@\ \ifcase\count@,
\or;\or.\or 
  \or\let\refspace@\relax
  \else,\fi
  \ifquotes@''\quotes@false\fi \penalty@ \refspace@
}
\def\transferpenalty@#1{\dimen@\lastkern\unkern
  \ifnum\lastpenalty=0\unpenalty\let\penalty@\relax
  \else\edef\penalty@{\penalty\the\lastpenalty\relax}\unpenalty\fi
  #1\penalty@\kern\dimen@}
\def\endref{\endrefitem@\lastref@true\endref@
  \par\endgroup \prevjour@false \previnbook@false }
\def\endref@{%
\iffirstref@
  \ifvoid\nobox@\ifvoid\keybox@\indent\fi
  \else\hbox to\refindentwd{\hss\unhbox\nobox@}\fi
  \ifvoid\keybox@
  \else\ifdim\wd\keybox@>\refindentwd
         \box\keybox@
       \else\hbox to\refindentwd{\unhbox\keybox@\hfil}\fi\fi
  \kern4sp\ppunbox@\bybox@
\fi 
  \ifvoid\paperbox@
  \else\prepunct@\unhbox\paperbox@
    \ifpaperquotes@\quotes@true\fi\fi
  \ppunbox@\paperinfobox@
  \ifvoid\jourbox@
    \ifprevjour@ \nocomma@\volbox@
      \nocomma@\issuebox@
      \ifvoid\yrbox@\else\changepunct@3\prepunct@(\unhbox\yrbox@
        \transferpenalty@)\fi
      \ppunbox@\pagesbox@
    \fi 
  \else \prepunct@\unhbox\jourbox@
    \nocomma@\volbox@
    \nocomma@\issuebox@
    \ifvoid\yrbox@\else\changepunct@3\prepunct@(\unhbox\yrbox@
      \transferpenalty@)\fi
    \ppunbox@\pagesbox@
  \fi 
  \ifbook@\prepunct@\unhbox\bookbox@ \ifbookquotes@\quotes@true\fi \fi
  \nocomma@\edsbox@
  \ppunbox@\bookinfobox@
  \ifbook@\ifvoid\volbox@\else\prepunct@ vol.~\unhbox\volbox@
  \fi\fi
  \ppunbox@\publbox@ \ppunbox@\publaddrbox@
  \ifbook@ \ppunbox@\yrbox@
    \ifvoid\pagesbox@
    \else\prepunct@\pp@\unhbox\pagesbox@\fi
  \else
    \ifprevinbook@ \ppunbox@\yrbox@
      \ifvoid\pagesbox@\else\prepunct@\pp@\unhbox\pagesbox@\fi
    \fi \fi
  \ppunbox@\finalinfobox@
  \iflastref@
    \ifvoid\langbox@.\ifquotes@''\fi
    \else\changepunct@2\prepunct@\unhbox\langbox@\fi
  \else
    \ifvoid\langbox@\changepunct@1%
    \else\changepunct@3\prepunct@\unhbox\langbox@
      \changepunct@1\fi
  \fi
}
\outer\def\enddocument{%
 \runaway@{proclaim}%
\ifmonograph@ 
\else
 \nobreak
 \thetranslator@
 \count@\z@ \loop\ifnum\count@<\addresscount@\advance\count@\@ne
 \csname address\number\count@\endcsname
 \csname email\number\count@\endcsname
 \repeat
\fi
 \vfill\supereject\end}

\def\headfont@{\headfonts}
\def\proclaimheadfont@{\bf}
\def\specialheadfont@{\bf}
\def\subheadfont@{\bf}
\def\demoheadfont@{\smc}

\newif\ifThisToToc \ThisToTocfalse
\newif\iftocloaded \tocloadedfalse

\def\C@L{\noexpand\Cal}\def\B@B{\noexpand\Bbb}\def\fR@K{\noexpand\frak}
\def\S@{\noexpand\S}\def\P@P{\noexpand\"}
\def\xpar{\\}

\def\writetoc#1{\iftocloaded\ifThisToToc\begingroup\def\totoc{}
  \def\Cal{\noexpand\C@L}\def\Bbb{\noexpand\B@B}
  \def\frak{\noexpand\fR@K}\def\goth{\frak}\def\S{\noexpand\S@}
  \def\"{\noexpand\P@P}
  \def\xpar{\par\penalty100000 }\def\idx##1{##1}\def\\{\xpar}
  \edef\next@{\write\toc{\noindent#1\leaderfill\noexpand\folio\par}}%
  \next@\endgroup\global\ThisToTocfalse\fi\fi}
\def\leaderfill{\leaders\hbox to 1em{\hss.\hss}\hfill}

\newif\ifindexloaded \indexloadedfalse
\def\idx#1{\ifindexloaded\begingroup\def\ign{}\def\it{}\def\/{}%
 \def\smc{}\def\bf{}\def\tt{}%
 \def\Cal{\noexpand\C@L}\def\Bbb{\noexpand\B@B}%
 \def\frak{\noexpand\fR@K}\def\goth{\frak}\def\S{\noexpand\S@}%
  \def\"{\noexpand\P@P}%
 {\edef\next@{\write\index{#1, \noexpand\folio}}\next@}%
 \endgroup\fi{#1}}
\def\ign#1{}

\def\totoc{\global\ThisToToctrue}

\outer\def\head#1\endhead{\par\penaltyandskip@{-200}\aboveheadskip
 {\headfont@\raggedcenter@\interlinepenalty\@M
 \ignorespaces#1\endgraf}\nobreak
 \vskip\belowheadskip
 \headmark{#1}\writetoc{#1}}

\outer\def\chaphead#1\endchaphead{\par\penaltyandskip@{-200}\aboveheadskip
 {\chapheadfonts\raggedcenter@\interlinepenalty\@M
 \ignorespaces#1\endgraf}\nobreak
 \vskip3\belowheadskip
 \headmark{#1}\writetoc{#1}}

\def\folio{{\foliofont@\ifnum\pageno<\z@ \romannumeral-\pageno
 \else\number\pageno \fi}}
\newtoks\leftheadtoks
\newtoks\rightheadtoks

\def\leftheadtext{\nofrills@{\relax}\lht@
  \DNii@##1{\leftheadtoks\expandafter{\lht@{##1}}%
    \mark{\the\leftheadtoks\noexpand\else\the\rightheadtoks}
    \ifsyntax@\setboxz@h{\def\\{\unskip\space\ignorespaces}%
        \headlinefont@##1}\fi}%
  \FN@\next@}
\def\rightheadtext{\nofrills@{\relax}\rht@
  \DNii@##1{\rightheadtoks\expandafter{\rht@{##1}}%
    \mark{\the\leftheadtoks\noexpand\else\the\rightheadtoks}%
    \ifsyntax@\setboxz@h{\def\\{\unskip\space\ignorespaces}%
        \headlinefont@##1}\fi}%
  \FN@\next@}
\def\NoRunningHeads{\global\runheads@false\global\let\headmark\eat@}

\newif\iffirstpage@     \firstpage@true
\newif\ifrunheads@      \runheads@true

\newdimen\fullhsize \fullhsize=\hsize
\newdimen\fullvsize \fullvsize=\vsize
\def\fullline{\hbox to\fullhsize}

\def\pagenumbers{\gdef\folio{\folio@}}

\let\norunningheads\NoRunningHeads
\def\userunningheads{\global\runheads@true}
\norunningheads

\headline={\def\chapter#1{\chapterno@. }%
  \def\\{\unskip\space\ignorespaces}\ifrunheads@\headlinefont@
    \ifodd\pageno\rightheadline \else\leftheadline\fi
   \else\hfil\fi\ifNoRunHeadline\global\NoRunHeadlinefalse\fi}
\let\folio@\folio
\def\foliofont@{\foliofont}
\def\foliofont{\eightrm}
\def\headlinefont@{\headlinefont}
\def\headlinefont{\eightpoint\smc}
\def\leftheadline{\rlap{\folio}\hfill
   \ifNoRunHeadline\else\iftrue\topmark\fi\fi \hfill}
\def\rightheadline{\hfill\ifNoRunHeadline
   \else \expandafter\fi
  \hfill \llap{\folio}}
\footline={{\eightpoint\bottremark}%
   \ifrunheads@\else\hfil{\let\foliofont\tenrm\folio}\fi\hfil}
\def\bottremark{}
 
\newif\ifNoRunHeadline      
\def\norunninghead{\global\NoRunHeadlinetrue}
\norunninghead

\output={\output@}
%
\newif\ifoffset\offsetfalse
\output={\output@}
\def\output@{%
 \ifoffset 
  \ifodd\count0\advance\hoffset by0.5truecm
   \else\advance\hoffset by-0.5truecm\fi\fi
 \shipout\vbox{%
  \makeheadline \pagebody \makefootline }%
 \advancepageno \ifnum\outputpenalty>-\@MM\else\dosupereject\fi}

\def\indexoutput#1{%
  \ifoffset 
   \ifodd\count0\advance\hoffset by0.5truecm
    \else\advance\hoffset by-0.5truecm\fi\fi
  \shipout\vbox{\makeheadline
  \vbox to\fullvsize{\boxmaxdepth\maxdepth%
  \ifvoid\topins\else\unvbox\topins\fi%
  #1 %
  \ifvoid\footins\else 
    \vskip\skip\footins
    \footnoterule
    \unvbox\footins\fi
  \ifr@ggedbottom \kern-\dimen@ \vfil \fi}%
  \baselineskip2pc
  \makefootline}%
 \global\advance\pageno\@ne
 \ifnum\outputpenalty>-\@MM\else\dosupereject\fi}
 
 \newbox\partialpage \newdimen\halfsize \halfsize=0.5\fullhsize
 \advance\halfsize by-0.5em

 \def\begindoublecolumns{\output={\indexoutput{\unvbox255}}%
   \begingroup \def\line{\fullline}
   \output={\global\setbox\partialpage=\vbox{\unvbox255\bigskip}}\eject
   \output={\doublecolumnout}\hsize=\halfsize \vsize=2\fullvsize}
 \def\enddoublecolumns{\output={\balancecolumns}\eject
  \endgroup \pagegoal=\fullvsize%
  \output={\output@}}
\def\doublecolumnout{\splittopskip=\topskip \splitmaxdepth=\maxdepth
  \dimen@=\fullvsize \advance\dimen@ by-\ht\partialpage
  \setbox0=\vsplit255 to \dimen@ \setbox2=\vsplit255 to \dimen@
  \indexoutput{\pagesofar} \unvbox255 \penalty\outputpenalty}
\def\pagesofar{\unvbox\partialpage
  \wd0=\hsize \wd2=\hsize \hbox to\fullhsize{\box0\hfil\box2}}
\def\balancecolumns{\setbox0=\vbox{\unvbox255} \dimen@=\ht0
  \advance\dimen@ by\topskip \advance\dimen@ by-\baselineskip
  \divide\dimen@ by2 \splittopskip=\topskip
  {\vbadness=10000 \loop \global\setbox3=\copy0
    \global\setbox1=\vsplit3 to\dimen@
    \ifdim\ht3>\dimen@ \global\advance\dimen@ by1pt \repeat}
  \setbox0=\vbox to\dimen@{\unvbox1} \setbox2=\vbox to\dimen@{\unvbox3}
  \pagesofar}

\tenpoint
\catcode`\@=\active

\def\smallheadings{\let\chapheadfonts\tenpoint\let\headfonts\tenpoint}

\tenpoint
\catcode`\@=\active

\def\LL{\leavevmode\setbox0=\hbox{L}\hbox to\wd0{\hss\char'40L}}
\def\al{\alpha}
\def\be{\beta}

\def\de{\delta}

\def\rh{\rho}

\def\om{\omega}


\def\today{\ifcase\month\or
 January\or February\or March\or April\or May\or June\or
 July\or August\or September\or October\or November\or December\fi
 \space\number\day, \number\year}

\def\({\left(}
\def\){\right)}
\def\[{\left[}
\def\]{\right]}

\def\Im{\operatorname{Im}}
\def\Re{\operatorname{Re}}

\def\3{\ss}
\catcode`\@=11
\def\dddot#1{\vbox{\ialign{##\crcr
      .\hskip-.5pt.\hskip-.5pt.\crcr\noalign{\kern1.5\p@\nointerlineskip}
      $\hfil\displaystyle{#1}\hfil$\crcr}}}

\newif\iftab@\tab@false
\newif\ifvtab@\vtab@false
\def\tab{\bgroup\tab@true\vtab@false\vst@bfalse\Strich@false%
   \def\\{\global\hline@@false%
     \ifhline@\global\hline@false\global\hline@@true\fi\cr}
   \edef\l@{\the\leftskip}\ialign\bgroup\hskip\l@##\hfil&&##\hfil\cr}
\def\endtab{\cr\egroup\egroup}
\def\vtab{\vtop\bgroup\vst@bfalse\vtab@true\tab@true\Strich@false%
   \bgroup\def\\{\cr}\ialign\bgroup&##\hfil\cr}
\def\endvtab{\cr\egroup\egroup\egroup}
\def\stab{\D@cke0.5pt\null 
 \bgroup\tab@true\vtab@false\vst@bfalse\Strich@true\Let@@\vspace@
 \normalbaselines\offinterlineskip
  \openup\spreadmlines@
 \edef\l@{\the\leftskip}\ialign
 \bgroup\hskip\l@##\hfil&&##\hfil\crcr}
\def\endstab{\crcr\egroup
 \egroup}
\newif\ifvst@b\vst@bfalse
\def\vstab{\D@cke0.5pt\null
 \vtop\bgroup\tab@true\vtab@false\vst@btrue\Strich@true\bgroup\Let@@\vspace@
 \normalbaselines\offinterlineskip
  \openup\spreadmlines@\bgroup}
\def\endvstab{\crcr\egroup\egroup
 \egroup\tab@false\Strich@false}

\newdimen\htstrut@
\htstrut@8.5\p@
\newdimen\htStrut@
\htStrut@12\p@
\newdimen\dpstrut@
\dpstrut@3.5\p@
\newdimen\dpStrut@
\dpStrut@3.5\p@
\def\openup{\afterassignment\@penup\dimen@=}
\def\@penup{\advance\lineskip\dimen@
  \advance\baselineskip\dimen@
  \advance\lineskiplimit\dimen@
  \divide\dimen@ by2
  \advance\htstrut@\dimen@
  \advance\htStrut@\dimen@
  \advance\dpstrut@\dimen@
  \advance\dpStrut@\dimen@}
\def\Let@@{\relax%
    \def\\{\global\hline@@false%
     \ifhline@\global\hline@false\global\hline@@true\fi\cr}%
    \iffalse}\fi}
\def\matrix{\null\,\vcenter\bgroup
 \tab@false\vtab@false\vst@bfalse\Strich@false\Let@@\vspace@
 \normalbaselines\openup\spreadmlines@\ialign
 \bgroup\hfil$\m@th##$\hfil&&\quad\hfil$\m@th##$\hfil\crcr
 \Mathstrut@\crcr\noalign{\kern-\baselineskip}}
\def\endmatrix{\crcr\Mathstrut@\crcr\noalign{\kern-\baselineskip}\egroup
 \egroup\,}
\def\smatrix{\D@cke0.5pt\null\,
 \vcenter\bgroup\tab@false\vtab@false\vst@bfalse\Strich@true\Let@@\vspace@
 \normalbaselines\offinterlineskip
  \openup\spreadmlines@\ialign
 \bgroup\hfil$\m@th##$\hfil&&\quad\hfil$\m@th##$\hfil\crcr}
\def\endsmatrix{\crcr\egroup
 \egroup\,\Strich@false}
\newdimen\D@cke
\def\Dicke#1{\global\D@cke#1}
\newtoks\tabs@\tabs@{&}
\newif\ifStrich@\Strich@false
\newif\iff@rst

\def\Stricherr@{\iftab@\ifvtab@\errmessage{\noexpand\s not allowed
     here. Use \noexpand\vstab!}%
  \else\errmessage{\noexpand\s not allowed here. Use \noexpand\stab!}%
  \fi\else\errmessage{\noexpand\s not allowed
     here. Use \noexpand\smatrix!}\fi}
\def\format{\ifvst@b\else\crcr\fi\egroup\iffalse{\fi\ifnum`}=0 \fi\format@}
\def\format@#1\\{\def\preamble@{#1}%
 \def\Str@chfehlt##1{\ifx##1\s\Stricherr@\fi\ifx##1\\\let\Next\relax%
   \else\let\Next\Str@chfehlt\fi\Next}%
 \def\c{\hfil\noexpand\ifhline@@\hbox{\vrule height\htStrut@%
   depth\dpstrut@ width\z@}\noexpand\fi%
   \ifStrich@\hbox{\vrule height\htstrut@ depth\dpstrut@ width\z@}%
   \fi\iftab@\else$\m@th\fi\the\hashtoks@\iftab@\else$\fi\hfil}%
 \def\r{\hfil\noexpand\ifhline@@\hbox{\vrule height\htStrut@%
   depth\dpstrut@ width\z@}\noexpand\fi%
   \ifStrich@\hbox{\vrule height\htstrut@ depth\dpstrut@ width\z@}%
   \fi\iftab@\else$\m@th\fi\the\hashtoks@\iftab@\else$\fi}%
 \def\l{\noexpand\ifhline@@\hbox{\vrule height\htStrut@%
   depth\dpstrut@ width\z@}\noexpand\fi%
   \ifStrich@\hbox{\vrule height\htstrut@ depth\dpstrut@ width\z@}%
   \fi\iftab@\else$\m@th\fi\the\hashtoks@\iftab@\else$\fi\hfil}%
 \def\s{\ifStrich@\ \the\tabs@\vrule width\D@cke\the\hashtoks@%
          \fi\the\tabs@\ }%
 \def\sa{\ifStrich@\vrule width\D@cke\the\hashtoks@%
            \the\tabs@\ %
            \fi}%
 \def\se{\ifStrich@\ \the\tabs@\vrule width\D@cke\the\hashtoks@\fi}%
 \def\cd{\hfil\noexpand\ifhline@@\hbox{\vrule height\htStrut@%
   depth\dpstrut@ width\z@}\noexpand\fi%
   \ifStrich@\hbox{\vrule height\htstrut@ depth\dpstrut@ width\z@}%
   \fi$\dsize\m@th\the\hashtoks@$\hfil}%
 \def\rd{\hfil\noexpand\ifhline@@\hbox{\vrule height\htStrut@%
   depth\dpstrut@ width\z@}\noexpand\fi%
   \ifStrich@\hbox{\vrule height\htstrut@ depth\dpstrut@ width\z@}%
   \fi$\dsize\m@th\the\hashtoks@$}%
 \def\ld{\noexpand\ifhline@@\hbox{\vrule height\htStrut@%
   depth\dpstrut@ width\z@}\noexpand\fi%
   \ifStrich@\hbox{\vrule height\htstrut@ depth\dpstrut@ width\z@}%
   \fi$\dsize\m@th\the\hashtoks@$\hfil}%
 \ifStrich@\else\Str@chfehlt#1\\\fi%
 \setbox\z@\hbox{\xdef\Preamble@{\preamble@}}\ifnum`{=0 \fi\iffalse}\fi
 \ialign\bgroup\span\Preamble@\crcr}
\newif\ifhline@\hline@false
\newif\ifhline@@\hline@@false
\def\hlinefor#1{\multispan@{\strip@#1 }\leaders\hrule height\D@cke\hfill%
    \global\hline@true\ignorespaces}
\def\Item "#1"{\par\noindent\hangindent2\parindent%
  \hangafter1\setbox0\hbox{\rm#1\enspace}\ifdim\wd0>2\parindent%
  \box0\else\hbox to 2\parindent{\rm#1\hfil}\fi\ignorespaces}
\def\ITEM #1"#2"{\par\noindent\hangafter1\hangindent#1%
  \setbox0\hbox{\rm#2\enspace}\ifdim\wd0>#1%
  \box0\else\hbox to 0pt{\rm#2\hss}\hskip#1\fi\ignorespaces}
\def\item"#1"{\par\noindent\hang%
  \setbox0=\hbox{\rm#1\enspace}\ifdim\wd0>\the\parindent%
  \box0\else\hbox to \parindent{\rm#1\hfil}\enspace\fi\ignorespaces}
\let\plainitem@\item
\catcode`\@=13

\hsize13cm
\vsize19cm
\newdimen\fullhsize
\newdimen\fullvsize
\newdimen\halfsize
\fullhsize13cm
\fullvsize19cm
\halfsize=0.5\fullhsize
\advance\halfsize by-0.5em

\magnification1200

\TagsOnRight

\def\ChriAA{1}
\def\ElouAA{2}
\def\IsmaAA{3}
\def\KratBN{4}
\def\KratCN{5}
\def\LascAZ{6}
\def\SzegAA{7}
\def\UvarAA{8}
\def\UvarAB{9}
\def\VienAE{10}

\def\AA{1.1}
\def\AB{1.2}
\def\AC{1.3}
\def\AE{1.4}
\def\AF{1.5}
\def\AG{1.6}
\def\BA{2.1}
\def\BB{2.2}
\def\AD{2.3}
\def\BBa{2.4}
\def\BBf{2.5}
\def\BBc{2.6}
\def\BBe{2.7}
\def\BBd{2.8}
\def\BBb{2.9}
\def\BC{3.1}
\def\BD{3.2}
\def\BE{3.3}

\def\CAa{4.1}
\def\CA{4.2}
\def\CB{4.3}
\def\CC{4.4}
\def\CD{4.5}
\def\CE{4.6}

\def\DA{5.1}

\def\ZA{7.1}
\def\ZB{7.2}
\def\ZC{7.3}
\def\ZD{7.4}
\def\ZE{7.5}
\def\ZF{7.6}
\def\ZG{7.7}
\def\ZH{7.8}
\def\ZI{7.9}
\def\ZIa{7.10}
\def\ZJ{7.11}
\def\ZK{7.12}
\def\ZL{7.13}
\def\ZM{7.14}
\def\ZN{7.15}
\def\ZO{7.16}
\def\ZP{7.17}
\def\ZQ{7.18}

%

\def\TA{1}
\def\TAa{2}
\def\TAb{3}
\def\TB{4}
\def\TC{5}
\def\TCa{6}
\def\TD{7}
\def\TE{8}
\def\TF{9}
\def\TGa{10}
\def\TG{11}
\def\TH{12}
\def\TI{13}
\def\TJ{14}
\def\TK{15}
\def\TL{16}

%


\def\fl#1{\lfloor#1\rfloor}
\def\cl#1{\lceil#1\rceil}

\topmatter 
\title A determinant identity for moments of orthogonal polynomials
that implies Uvarov's formula for the orthogonal polynomials 
of rationally related densities
\endtitle 
\author C.~Krattenthaler
\endauthor 
\affil 
Fakult\"at f\"ur Mathematik, Universit\"at Wien,\\
Oskar-Morgenstern-Platz~1, A-1090 Vienna, Austria.\\
WWW: \tt http://www.mat.univie.ac.at/\~{}kratt
\endaffil
\address Fakult\"at f\"ur Mathematik, Universit\"at Wien,
Oskar-Morgenstern-Platz~1, A-1090 Vienna, Austria.\newline
http://www.mat.univie.ac.at/\~{}kratt
\endaddress

\thanks Research partially supported by the Austrian
Science Foundation FWF (grant S50-N15)
in the framework of the Special Research Program
``Algorithmic and Enumerative Combinatorics"%
\endthanks

\subjclass Primary 33C45;
 Secondary 05A10 05A19 11C20 15A15\linebreak 42C05
\endsubjclass
\keywords Hankel determinants, moments of orthogonal polynomials,
Dodgson con\-den\-sa\-tion, Chebyshev polynomials, Catalan numbers
\endkeywords
\abstract 
Let $p_n(x)$, $n=0,1,\dots$, be the orthogonal polynomials with
respect to a given density $d\mu(x)$. Furthermore, let $d\nu(x)$
be a density which arises from $d\mu(x)$ by multiplication by a
rational function in~$x$. We prove a formula that expresses the
Hankel determinants of moments of $d\nu(x)$ in terms of a
determinant involving the orthogonal polynomials $p_n(x)$
and associated functions $q_n(x)=\int p_n(u) \,d\mu(u)/(x-u)$.
Uvarov's formula for the orthogonal polynomials with respect to
$d\nu(x)$ is a corollary of our theorem. Our result generalises
a Hankel determinant formula for the case where the rational
function is a polynomial that existed somehow hidden in the folklore
of the theory of orthogonal polynomials but has been stated explicitly only
relatively recently (see [{\tt ar$\chi$iv:2101.04225}]).
Our theorem can be interpreted in a two-fold way: analytically
or in the sense of formal series.
We apply our theorem to derive several curious Hankel determinant
evaluations.
\endabstract
\endtopmatter
\document

\subhead 1. Introduction\endsubhead
Recently, in \cite{\KratCN} this author discovered a formula that
expresses the Hankel determinant of linear combinations of 
moments of orthogonal polynomials
in terms of a determinant involving these orthogonal polynomials.
A literature search revealed that this formula existed in a hidden
form behind a theorem (cf\. \cite{\SzegAA,
Theorem~2.5} or \cite{\IsmaAA, Theorem~2.7.1})
that is commonly attributed to Christoffel
\cite{\ChriAA} (although he had only proved it in a very special case);
only recently it had been stated explicitly, by Lascoux
in \cite{\LascAZ, Prop.~8.4.1} (although incorrectly) and
by Elouafi \cite{\ElouAA, Theorem~1} (however with an incomplete
proof). Three fundamentally different
proofs are given in \cite{\KratCN}: one due to
this author, one following Lascoux's arguments, and one completing
Elouafi's arguments.

The purpose of this article is to present and prove a generalisation
of the aforementioned formula that is inspired by Uvarov's formula
\cite{\UvarAA, \UvarAB} for the orthogonal polynomials with respect
to a density that is related to another given density by the
multiplication by a rational function, see Theorem~\TA\ below.

\medskip
Let $\big(p_n(x)\big)_{n\ge0}$ be a sequence of monic polynomials
over a field~$K$ of characteristic zero\footnote{For the analyst, 
(usually) this field is the field of 
real numbers, and a further restriction is that the linear
functional~$L$ is defined by a measure with non-negative density. 
However, the formulae
in this paper do not need these restrictions and 
are valid in this wider context of ``formal orthogonality".} with
$\deg p_n(x)=n$, and assume that they are orthogonal with respect to the linear
functional~$L$, i.e., they satisfy $L(p_m(x)p_n(x))=\om_n\de_{m,n}$ with
$\om_n\ne0$ for all~$n$,
where $\de_{m,n}$ is the Kronecker delta.
Furthermore, we write $\mu_n$ for the $n$-th moment
$L(x^n)$ of the functional~$L$. For convenience (and in order to keep
the usual analytic meaning in mind), we shall write
$\int f(u)\,d\mu(u)$ instead of $L(f(x))$. This can be either read in
a purely formal way, or the analyst may think of it as a concrete
integral with respect to the measure given by the density $d\mu(u)$.

For the statement of our theorem, we need the ``functions"
$$
q_n(y)=\int \frac {p_n(u)} {y-u}\,d\mu(u).
\tag\AA
$$
These can be understood in the ``ordinary" analytic sense if
$\mu(u)$ is a concrete measure, or, alternatively, these can be
understood in the sense of formal power series in $1/y$,
see Lemma~\TCa, Equation~(\BBe
).

\medskip
Here is the main result of this article.

\proclaim{Theorem \TA}
Let $k$, $m$ and $n$ be non-negative integers
and $x_1,x_2,\dots,x_m$ and $y_1,y_2,\mathbreak\dots,y_k$ be variables.
Then, with the above notations, for $n\ge k$ we have
$$
\multline
\frac {\det\limits_{0\le i,j\le n-1}\left(
\dsize\int u^{i+j}\dfrac {\prod _{\ell=1} ^{m}(u-x_\ell)}
{\prod _{\ell=1} ^{k}(u-y_\ell)}\,d\mu(u)\right)}
{\det\limits_{0\le i,j\le n-k-1}\left(\mu_{i+j}\right)}
\\
=(-1)^{n(m-k)+km}
\frac {\displaystyle\det M_{k,m,n}(x_1,\dots,x_m,y_1,\dots,y_k)}
{
\bigg(\prod\limits _{1\le i<j\le m} ^{}(x_j-x_i)\bigg)
\bigg(\prod\limits _{1\le i<j\le k} ^{}(y_i-y_j)\bigg)
},
\endmultline
\tag\AB
$$
where
$$
M_{k,m,n}(x_1,\dots,x_m,y_1,\dots,y_k)
=
\pmatrix 
p_{n-k}(x_1)&p_{n-k+1}(x_1)&\dots&p_{n+m-1}(x_1)\\
\hdotsfor4\\
p_{n-k}(x_m)&p_{n-k+1}(x_m)&\dots&p_{n+m-1}(x_m)\\
q_{n-k}(y_1)&q_{n-k+1}(y_1)&\dots&q_{n+m-1}(y_1)\\
\hdotsfor4\\
q_{n-k}(y_k)&q_{n-k+1}(y_k)&\dots&q_{n+m-1}(y_k)
\endpmatrix.
$$
If $n< k$, then
$$
\multline
{\det\limits_{0\le i,j\le n-1}\left(
\dsize\int u^{i+j}\dfrac {\prod _{\ell=1} ^{m}(u-x_\ell)}
{\prod _{\ell=1} ^{k}(u-y_\ell)}\,d\mu(u)\right)}
\\
=(-1)^{n(m-k)+km}
\frac {\displaystyle\det N_{k,m,n}(x_1,\dots,x_m,y_1,\dots,y_k)}
{
\bigg(\prod\limits _{1\le i<j\le m} ^{}(x_j-x_i)\bigg)
\bigg(\prod\limits _{1\le i<j\le k} ^{}(y_i-y_j)\bigg)
},
\endmultline
\tag\AC
$$
where
$$
\multline
N_{k,m,n}(x_1,\dots,x_m,y_1,\dots,y_k)
\\
=
\pmatrix 
0&\dots&0&0&0&p_{0}(x_1)&p_{1}(x_1)&\dots&p_{n+m-1}(x_1)\\
\hdotsfor9\\
0&\dots&0&0&0&p_{0}(x_m)&p_{1}(x_m)&\dots&p_{n+m-1}(x_m)\\
y_1^{k-n-1}&\dots&y_1^2&y_1&1&q_{0}(y_1)&q_{1}(y_1)&\dots&q_{n+m-1}(y_1)\\
\hdotsfor9\\
y_k^{k-n-1}&\dots&y_k^2&y_k&1&q_{0}(y_k)&q_{1}(y_k)&\dots&q_{n+m-1}(y_k)
\endpmatrix.
\endmultline
$$
Here, determinants of empty matrices and empty products are understood
to equal~$1$.
\endproclaim

\remark{Remarks}
(1) The numerator determinant on the left-hand sides of (\AB) and~(\AC) 
is the Hankel determinant $\det_{0\le i,j\le n-1}(\rh_{i+j})$, where
the $\rh_s$'s are the moments of the linear functional
$$
p(x)\mapsto \int p(u)\dfrac {\prod _{\ell=1} ^{m}(u-x_\ell)}
{\prod _{\ell=1} ^{k}(u-y_\ell)}\,d\mu(u).
$$

\medskip
(2) The theory of orthogonal polynomials guarantees that in our setting
(namely due to the condition $\om_n\ne0$ in the orthogonality)
the Hankel determinant of moments in the denominator on the left-hand
side of~(\AB) is non-zero.

\medskip
(3) The main theorem in \cite{\KratCN} is equivalent to the special
case of Theorem~ \TA\ where $k=0$. It is worth stating this special
case separately.

\proclaim{Corollary \TAa}
Let $m$ and $n$ be non-negative integers
and $x_1,x_2,\dots,x_m$ be variables.
Then we have
$$
\frac {\det\limits_{0\le i,j\le n-1}\left(
\dsize\int  {u^{i+j}}
{\prod _{\ell=1} ^{m}(u-x_\ell)}\,d\mu(u)\right)}
{\det\limits_{0\le i,j\le n-1}\left(\mu_{i+j}\right)}
=(-1)^{nm}
\frac {\displaystyle\det_{1\le i,j\le m}\big(p_{n+j-1}(x_i)\big)}
{
\prod\limits _{1\le i<j\le m} ^{}(x_j-x_i)
}.
\tag\AE
$$
\endproclaim

\medskip
(4) Similarly, it is worth stating the special case of Theorem~\TA\
where $m=0$ separately.

\proclaim{Corollary \TAb}
Let $k$ and $n$ be non-negative integers
and $y_1,y_2,\dots,y_k$ be variables.
Then, for $n\ge k$ we have
$$
\frac {\det\limits_{0\le i,j\le n-1}\left(
\dsize\int \dfrac {u^{i+j}\,d\mu(u)}
{\prod _{\ell=1} ^{k}(u-y_\ell)}\right)}
{\det\limits_{0\le i,j\le n-k-1}\left(\mu_{i+j}\right)}
=(-1)^{nk}
\frac {\displaystyle\det_{1\le i,j\le k}\big(q_{n-k+j-1}(y_i)\big) }
{
\prod\limits _{1\le i<j\le k} ^{}(y_i-y_j)
}.
\tag\AF
$$
If $n< k$, then
$$
{\det\limits_{0\le i,j\le n-1}\left(
\dsize\int \dfrac {u^{i+j}\,d\mu(u)}
{\prod _{\ell=1} ^{k}(u-y_\ell)}\right)}
=(-1)^{nk}
\frac {\displaystyle\det N_{k,n}(y_1,\dots,y_k)}
{
\prod\limits _{1\le i<j\le k} ^{}(y_i-y_j)
},
\tag\AG
$$
where
$$
N_{k,n}(y_1,\dots,y_k)
=
\pmatrix 
y_1^{k-n-1}&\dots&y_1^2&y_1&1&q_{0}(y_1)&q_{1}(y_1)&\dots&q_{n-1}(y_1)\\
\hdotsfor9\\
y_k^{k-n-1}&\dots&y_k^2&y_k&1&q_{0}(y_k)&q_{1}(y_k)&\dots&q_{n-1}(y_k)
\endpmatrix.
$$
\endproclaim

\medskip
(5) Let $d\mu(u)$ be a given density.
We will explain in Section~6 how Theorem~\TA\ is related to
Uvarov's formula \cite{\UvarAA, \UvarAB} for the orthogonal
polynomials with respect to the linear functional defined by the
density
$$
\dfrac {\prod _{\ell=2} ^{m}(u-x_\ell)}
{\prod _{\ell=1} ^{k}(u-y_\ell)}\,d\mu(u).
$$

\medskip
(6) In \cite{\KratCN}, three proofs of the special case of
Theorem~\TA\ where $k=0$ --- that is, of Corollary~\TAa\ ---
are given, one using the method of condensation, one using classical
results from the theory of orthogonal polynomials, and one using a 
vanishing argument. It is interesting to note that neither the
second nor the third proof seem to extend to a proof of Theorem~\TA,
only the condensation argument does. This is indeed the argument
that we apply here in Section~4.
\endremark

The rest of this article is organised as follows. In the next section,
we review some classical facts from the theory of orthogonal
polynomials that will be relevant for the proof of Theorem~\TA.
Furthermore, in Lemma~\TCa, we provide some information on the
``functions" $q_n(y)$ that are so central in Theorem~\TA.
Our proof of Theorem~\TA\ requires several determinant identities.
These are presented in Section~3, among which Jacobi's condensation
formula (see Lemma~\TE). The actual proof of Theorem~\TA\ is then
the subject of Section~4. In Proposition \TI\ in Section~5 it is
explained how to ``read" Theorem~\TA\ in cases where two (or more)
of the $x_i$'s or the $y_i$'s are equal. We come back to the title
of this article in Section~6 by outlining how Theorem~\TA\ relates 
to Uvarov's formula~\cite{\UvarAA, \UvarAB}. Finally, by
considering the special case of Theorem~\TA\ where
$d\mu(u)=\sqrt{1-u^2}\,du$ --- corresponding to Chebyshev
polynomials of the second kind ---, 
we derive evaluations of several curious Hankel determinants
in Section~7.

\subhead 2. Preliminaries on orthogonal polynomials\endsubhead
In this section we survey classical facts about orthogonal
polynomials that we shall need in the sequel. We also prove
some properties of the ``functions" $q_n(y)$.

\medskip
As in the introduction, let $\big(p_n(x)\big)_{n\ge0}$ be a sequence
of monic polynomials over a field~$K$ of characteristic zero that is
orthogonal with respect to the linear function~$L$ given by
$$
L:p(x)\mapsto \int p(u)\,d\mu(u),
$$
where this may be read formally, or --- when we want to interpret this
analytically --- where $d\mu(u)$ is some given density. Such a sequence
of orthogonal polynomials exists if and only if all Hankel
determinants $\det_{0\le i,j\le n-1}(\mu_{i+j})$ of moments
$\mu_s=\int u^s\,d\mu(u)$ do not vanish. For explicit formulae
for the orthogonal polynomials $p_n(x)$ in terms of the moments
see Lemmas~\TB\ and~\TC\ below.

By Favard's theorem (see e.g\. \cite{\KratBN,
Theorems~11--13}), the sequence $\big(p_n(x)\big)_{n\ge0}$ is
orthogonal if and only if it satisfies a
{\it three-term recurrence}
$$p_n(x)=(x-s_{n-1})p_{n-1}(x)-t_{n-2}p_{n-2}(x),
\quad \text{for }n\ge1,
\tag\BA
$$
with initial values
$p_{-1} (x)= 0$ and $p_0 (x)=1$, 
for some sequences $(s_n)_{n\ge0}$ and $(t_n)_{n\ge0}$
of elements of~$K$ with $t_n\ne0$ for all $n$.

The $t_n$'s are connected with the Hankel determinants of moments
by the formula (see e.g\. \cite{\VienAE, Ch.~IV, Cor.~6}) 
$$
\det_{0\le i,j\le n-1}\left(\mu_{i+j}\right)
=\prod _{i=0} ^{n-1}t_i^{n-i-1}.
\tag\BB
$$
Conversely, we may use this formula to express $t_{n-1}$ in terms of
the Hankel determinants of moments.
For convenience, we introduce the abbreviation
$$
H(n):=\det\limits_{0\le i,j\le n-1}\left(\mu_{i+j}\right).
\tag\AD
$$
Then, from (\BB) we infer
$$
t_{n-1}=\frac {H(n+1)/H(n)} {H(n)/H(n-1)}.
\tag\BBa$$

Next we quote two formulae that express orthogonal polynomials
in terms of their associated moments.
The first can be found in \cite{\SzegAA, p.~27, Eq.~(2.2.6)},
and the second in \cite{\SzegAA, p.~27, Eq.~(2.2.9)}).

\proclaim{Lemma \TB}
Let $M$ be a linear functional on polynomials in $x$ with
moments~$\nu_n$, $n=0,1,\dots$, such that all Hankel determinants 
$\det_{0\le i,j\le n}(\nu_{i+j})$, $n=0,1,\dots$, 
are non-zero. 
Then the determinants
$$
\frac {1} {\det\limits_{0\le i,j\le n-1}\left(
\nu_{i+j}\right)}\det\pmatrix 
\nu_0&\nu_1&\nu_2&\dots&\nu_n\\
\nu_1&\nu_2&\nu_3&\dots&\nu_{n+1}\\
\hdotsfor5\\
\nu_{n-1}&\nu_n&\nu_{n+1}&\dots&\nu_{2n-1}\\
1&x&x^2&\dots&x^n
\endpmatrix
$$
are a sequence of monic orthogonal polynomials with respect to~$M$.
\endproclaim

\proclaim{Lemma \TC}
Let $M$ be a linear functional on polynomials in $x$ with
moments~$\nu_n$, $n=0,1,\dots$, such that all Hankel determinants 
$\det_{0\le i,j\le n}(\nu_{i+j})$, $n=0,1,\dots$, 
are non-zero. 
Then the determinants
$$
\det_{0\le i,j\le n-1}\left(\nu_{i+j+1}-\nu_{i+j}x\right)
\tag\BBf$$
are a sequence of orthogonal polynomials with respect to~$M$.
\endproclaim

We now turn to the ``functions" $q_n(y)$ defined in~(\AA).
They satisfy the same three-term recurrence relation as the original
orthogonal polynomials $p_n(x)$. Moreover, it follows from
Lemma~\TB\ that $q_n(y)$ can be seen as (formal) power
series in $1/y$ of degree~$-n-1$.\footnote{By ``degree" of a formal
power series $f(y)$ in $1/y$ we mean the maximal exponent~$e$ such that $y^e$
appears in $f(y)$ with non-zero coefficient.
We warn the reader that
the recurrence (\BBc) is deceiving in regard of the degree of
$q_n(y)$: it seems to suggest that --- similar to the situation for
the underlying orthogonal polynomials $p_n(x)$ --- the degree would
rise by~1 when going from $q_{n-1}(y)$ to $q_n(y)$. However, as 
(\BBe) shows, on the contrary the degree {\it drops} by~1, caused
by a cancellation
of leading terms in~(\BBc).} 

\proclaim{Lemma \TCa}
We have
$$q_n(y)=(y-s_{n-1})q_{n-1}(y)-t_{n-2}q_{n-2}(y),
\quad \text{for }n\ge2,
\tag\BBc
$$
with initial values
$q_0 (y)=\int\frac {d\mu(u)} {y-u}$ and
$q_1 (y)=(y-s_0)\int\frac {d\mu(u)} {y-u}-\mu_0$, 
with the sequences $(s_n)_{n\ge0}$ and $(t_n)_{n\ge0}$
of elements of~$K$ that feature in the three-term recurrence {\rm(\BA)} for
the underlying orthogonal polynomials $p_n(x)$.

Moreover,
for all non-negative integers $n$, as a formal power series in~$1/y$
the ``function" $q_n(y)$ starts as
$$
q_n(y)=\frac {H(n+1)} {H(n)}
y^{-n-1}+
O\left(y^{-n-2}\right),
\tag\BBe$$
where $H(n)$ and $H(n+1)$ is the short notation introduced in {\rm(\AD)}.
\endproclaim

\demo{Proof}
Let $n\ge2$.
From (\BA), we get
$$\int \frac {p_n(u)} {y-u}\,d\mu(u)
=\int \frac {(u-s_{n-1})p_{n-1}(u)} {y-u}\,d\mu(u)
-t_{n-2}\int \frac{p_{n-2}(u)} {y-u}\,d\mu(u).
\tag\BBd$$
The first term on the right-hand side can be simplified
as follows:
$$\align
\int \frac {(u-s_{n-1})p_{n-1}(u)} {y-u}\,d\mu(u)
&=
\int \frac {(y-s_{n-1})p_{n-1}(u)} {y-u}\,d\mu(u)
-\int \frac {(y-u)p_{n-1}(u)} {y-u}\,d\mu(u)\\
&=
(y-s_{n-1})q_{n-1}(y).
\endalign$$
If this is used in (\BBd) together with the definition (\AA) of
$q_n(y)$, the recurrence~(\BBc) results immediately.
The initial values for $q_0(y)$ and $q_1(y)$ are straightforward to
derive from $p_0(x)=1$ and $p_1(x)=x-s_0$.

\medskip
In order to show the second assertion, we note that,
by definition of $q_n(y)$, we have
$$
q_n(y)=\int \frac {p_n(u)} {y-u}\,d\mu(u)
=\sum_{i=0}^\infty\int p_n(u)u^iy^{-i-1}\,d\mu(u).
$$
Because of the orthogonality of $p_n(u)$ with respect to the density
$d\mu(u)$, all terms of the above sum with $i<n$ vanish. Thus,
$$
q_n(y)
=\sum_{i=n}^\infty\int p_n(u)u^iy^{-i-1}\,d\mu(u)
=y^{-n-1}\int p_n(u)u^n\,d\mu(u)+O\left(y^{-n-2}\right).
\tag\BBb
$$
Now we use the formula of Lemma~\TB\ with $\nu_s=\mu_s$ for
all~$s$, to obtain
$$\align 
\int p_n(u)u^n\,d\mu(u)
&=\frac {1} {H(n)}\int 
\det\pmatrix 
\mu_0&\mu_1&\mu_2&\dots&\mu_n\\
\mu_1&\mu_2&\mu_3&\dots&\mu_{n+1}\\
\hdotsfor5\\
\mu_{n-1}&\mu_n&\mu_{n+1}&\dots&\mu_{2n-1}\\
1&u&u^2&\dots&u^n
\endpmatrix
u^n\,d\mu(u)\\
&=\frac {1} {H(n)}
\pmatrix 
\mu_0&\mu_1&\mu_2&\dots&\mu_n\\
\mu_1&\mu_2&\mu_3&\dots&\mu_{n+1}\\
\hdotsfor5\\
\mu_{n-1}&\mu_n&\mu_{n+1}&\dots&\mu_{2n-1}\\
\mu_n&\mu_{n+1}&\mu_{n+2}&\dots&\mu_{2n}
\endpmatrix
=\frac {H(n+1)} {H(n)}.
\endalign$$
If this is substituted back in (\BBb), the assertion of the lemma
follows immediately.\quad \quad \qed
\enddemo

\subhead 3. Auxiliary determinant identities\endsubhead
The purpose of this section is to collect three determinant formulae
that will turn out to be crucial in our proof of Theorem~\TA.

The proof method of our proof in Section~4 is
the {\it method of condensation} (frequently referred to as ``Dodgson
condensation"; see \cite{\KratBN, Sec.~2.3}). 
This method provides inductive proofs
that are based on a determinant identity due to Jacobi, which
we recall in the following proposition.

\proclaim{Proposition \TD}
Let $A$ be an $N\times N$ matrix. Denote the submatrix of $A$ in which
rows $i_1,i_2,\dots,i_k$ and columns $j_1,j_2,\dots,j_k$ are 
omitted by $A_{i_1,i_2,\dots,i_k}^{j_1,j_2,\dots,j_k}$. Then we have
$$
\det A\cdot \det A_{i_1,i_2}^{j_1,j_2}=\det A_{i_1}^{j_1}\cdot 
\det A_{i_2}^{j_2}-
\det A_{i_1}^{j_2}\cdot \det A_{i_2}^{j_1}
\tag\BC
$$
for all integers $i_1,i_2,j_1,j_2$ with $1\le i_1<i_2\le N$
and $1\le j_1<j_2\le N$.
\endproclaim

We need two further determinantal formulae, which involve
Hankel determinants of linear combinations of sequence elements.

The following is \cite{\KratCN, Lemma~3 with $n$ replaced by $n-1$}.

\proclaim{Lemma \TE}
Let $(c_n)_{n\ge0}$ be a given sequence, and $\al$ and $\be$ be
variables. Then, for all positive integers $n$, we have
$$\multline 
(\be-\al)\det_{0\le i,j\le n-2}\big(\al\be
c_{i+j}+(\al+\be)c_{i+j+1}+c_{i+j+2}\big)
\det_{0\le i,j\le n-1}\big(c_{i+j}\big)\\
=
\det_{0\le i,j\le n-2}\big(\al c_{i+j}+c_{i+j+1}\big)
\det_{0\le i,j\le n-1}\big(\be c_{i+j}+c_{i+j+1}\big)\\
-
\det_{0\le i,j\le n-2}\big(\be c_{i+j}+c_{i+j+1}\big)
\det_{0\le i,j\le n-1}\big(\al c_{i+j}+c_{i+j+1}\big).
\endmultline
\tag\BD$$
\endproclaim

We require another, similarly looking determinant identity.

\proclaim{Lemma \TF}
Let $(c_n)_{n\ge0}$ be a given sequence, and $\al$ and $\be$ be
variables. Then, for all positive integers $n$, we have
$$\multline 
\det_{0\le i,j\le n-1}\big(\al
c_{i+j}+c_{i+j+1}\big)
\det_{0\le i,j\le n-1}\big(\be
c_{i+j}+c_{i+j+1}\big)\\
=
-\det_{0\le i,j\le n}\big(c_{i+j}\big)
\det_{0\le i,j\le n-2}\big(\al\be
c_{i+j}+(\al+\be)c_{i+j+1}+c_{i+j+2}\big)\\
+\det_{0\le i,j\le n-1}\big(c_{i+j}\big)
\det_{0\le i,j\le n-1}\big(\al\be
c_{i+j}+(\al+\be)c_{i+j+1}+c_{i+j+2}\big).
\endmultline
\tag\BE$$
\endproclaim

\demo{Proof}
In principle, it should be possible to prove this directly in the
style of the proof of Lemma~\TE\ given in \cite{\KratCN}, that is, by
extracting the coefficient of $\al^s\be^t$ on both sides of~(\BE),
and by then reducing everything to some known determinant
identity. Embarrassingly, I failed to carry this through.
Therefore, instead I take recourse to a dirty trick that is based
on the fact that Theorem~\TA\ had been proved earlier for the
special case in which $k=0$ (cf\. \cite{\KratCN, Theorem~1} and 
Proposition~\TG\ below). Namely, we interpret the sequence
$(c_{n})_{n\ge0}$ as the sequence of moments of some linear
functional, and we let $\big(r_n(x)\big)_{n\ge0}$ be the
corresponding sequence of monic orthogonal polynomials with respect
to that functional. As was said earlier, for the orthogonal
polynomials to exist we must assume that all Hankel determinants
$\det_{0\le i,j\le n-1}(c_{i+j})$ are non-zero, which we do for
the moment. Since both sides of~(\BE) are polynomials in the $c_i$'s,
$\al$ and~$\be$, this restriction can be removed in the end.

Now, in the above setting, by Theorem~ \TA\ with $k=0$, $m=1$,
and $x_1=-\al$, we have
$$
\det_{0\le i,j\le n-1}\big(\al
c_{i+j}+c_{i+j+1}\big)=
H(n)r_n(-\al),
$$
and, by Theorem~ \TA\ with $k=0$, $m=2$, $x_1=-\al$
and $x_2=-\be$, we have
$$
\det_{0\le i,j\le n-1}\big(\al\be
c_{i+j}+(\al+\be)c_{i+j+1}+c_{i+j+2}\big)
=
H(n)\frac {r_n(-\al)r_{n+1}(-\be)-r_{n+1}(-\al)r_n(-\be)} {\al-\be}.
$$
Consequently, Equation~(\BE) turns out to be equivalent with
$$\multline 
\big(\tilde H(n)\big)^2r_n(-\al)r_n(-\be)
=
\tilde H(n+1)\tilde H(n-1)\frac {r_{n-1}(-\al)r_{n}(-\be)
-r_{n}(-\al)r_{n-1}(-\be)}
{\al-\be}\\
+\big(\tilde H(n)\big)^2
\frac {r_n(-\al)r_{n+1}(-\be)-r_{n+1}(-\al)r_n(-\be)} {\al-\be}.
\endmultline$$
Here we used the short notation $\tilde H(n)=\det_{0\le i,j\le n-1}(c_{i+j})$.
After cancellation of common factors using~(\BB) (with $\tilde H(n)$
in place of $H(n)$ and the appropriate sequence $(\tilde t_n)_{n\ge0}$
instead of $(t_n)_{n\ge0}$), this becomes
$$\multline 
(\al-\be)
r_n(-\al)r_n(-\be)
=
\tilde t_{n-1}\big( {r_{n-1}(-\al)r_{n}(-\be)-r_{n}(-\al)r_{n-1}(-\be)}\big)
\\
+
\big( {r_n(-\al)r_{n+1}(-\be)-r_{n+1}(-\al)r_n(-\be)}\big).
\endmultline$$
Indeed, this is trivially true because of~(\BA) with $n$ replaced by
$n+1$, $p_n(x)$ replaced by $r_n(x)$, $s_n$ replaced by $\tilde s_n$,
$t_n$ replaced by $\tilde t_n$,
and $x=-\be$, respectively $x=-\al$.\quad \quad \qed
\enddemo

\subhead 4. Proof of Theorem \TA\endsubhead
Here we prove Theorem~\TA. As announced, we shall apply the
method of condensation (see Proposition~\TD) in order to set up
an inductive proof. Propositions~\TG\ and~\TH\ will provide the
start of the induction. In order to ``keep the induction running",
we need a non-vanishing result that we present first.

\proclaim{Lemma \TGa}
Under the assumption that the Hankel determinants $H(n)$ in~{\rm(\AD)}
do not vanish, 
the numerator determinant on the left-hand sides of\/
{\rm(\AB)} and\/ {\rm(\AC)},
$$
\det\limits_{0\le i,j\le n-1}\left(
\dsize\int u^{i+j}\dfrac {\prod _{\ell=1} ^{m}(u-x_\ell)}
{\prod _{\ell=1} ^{k}(u-y_\ell)}\,d\mu(u)\right),
\tag\CAa$$ 
does not vanish identically.
\endproclaim

\demo{Proof}
Adopting the viewpoint of formal series in the variables
$x_1,x_2,\dots,x_m$ and $1/y_1,1/y_2,\dots,1/y_k$, the highest degree
term in~(\CAa) is
$$\align
\det\limits_{0\le i,j\le n-1}\left(
(-1)^{m-k}\mu_{i+j}\dfrac {\prod _{\ell=1} ^{m}x_\ell}
{\prod _{\ell=1} ^{k}y_\ell}\right)
&=(-1)^{n(m-k)}\frac {\prod _{j=1} ^{m}x_j^n} {\prod _{j=1} ^{k}y_j^n}
\det\limits_{0\le i,j\le n-1}\left(
\mu_{i+j}\right)\\
&=(-1)^{n(m-k)}H(n)
\frac {\prod _{j=1} ^{m}x_j^n} {\prod _{j=1} ^{k}y_j^n},
\endalign
$$
which is non-zero by our assumption of non-vanishing of $H(n)$. 

If we adopt the analytic point of view, then one would multiply
the determinant (\CAa) by $ {\prod _{j=1} ^{m}x_j^{-n}} \big/{\prod
_{j=1} ^{k}y_j^{-n}}$ and then compute the limit\footnote{This limit
is unproblematic if $d\mu(u)$ is a measure with finite
support.} as
$x_i\to\infty$ and $y_j\to \infty$ for all~$i$ and~$j$.
The result would be the determinant
$$
\det\limits_{0\le i,j\le n-1}\left(
(-1)^{m-k}\mu_{i+j}
\right)=(-1)^{n(m-k)}H(n),
$$
with the same conclusion.\quad \quad \qed
\enddemo

\proclaim{Proposition \TG}
Theorem {\rm\TA} holds for $k=0$.
\endproclaim

This is the main theorem (namely Theorem~1) in \cite{\KratCN}, for
which three fundamentally different proofs are provided there.

\proclaim{Proposition \TH}
Theorem {\rm\TA} holds for $m=0$.
\endproclaim

\demo{Proof}
We prove the claim by induction on $k$. 

\medskip
For the start of the induction we need the validity of
(\AB) for $k=m=0$ --- this is obvious since both left-hand and
right-hand side equal~1 in that case ---, of (\AC) for $m=0$ and $k=1$
--- also this is easy to see since the only case to consider is $n=0$,
which makes the left-hand side reduce to~1, and also the right-hand
side, due to the evaluation of the Vandermonde determinant ---, 
and of (\AB) for $m=0$ and $k=1$. The latter needs an argument,
which we provide next.

For $m=0$ and $k=1$, the numerator of the left-hand side of~(\AB)
reads 
$$
\det\limits_{0\le i,j\le n-1}\left(
\dsize\int u^{i+j}\dfrac {d\mu(u)}
{u-y_1}\right).
\tag\CA$$
We have
$$
\int u^{i+j}\dfrac {d\mu(u)} {u-y_1}
=
\int u^{i+j-1}\dfrac {u-y_1+y_1} {u-y_1}\,d\mu(u)
=\mu_{i+j-1}+y_1
\int u^{i+j-1}\dfrac {d\mu(u)} {u-y_1}.
$$
We use this relation in the last row, that is, for $i=n-1$.
Subsequently, we subtract the $(n-2)$-nd row multiplied by~$y_1$
from the last row (the $(n-1)$-st row). Thereby the entry in the
$j$-th column of the last row becomes $\mu_{n+j-2}$.
We repeat this operation with the $(n-2)$-nd and the $(n-3)$-rd row,
with the $(n-3)$-rd and the $(n-4)$-th row, \dots, and with the
first row and the $0$-th row. As a result, we have converted the
determinant in~(\CA) into the determinant
$$\align
\det&\pmatrix 
\dsize\int u^{j}\dfrac {d\mu(u)}
{u-y_1}&\text{for }i=0,\ 0\le j\le n-1\hfill\\
\mu_{i+j-1}&\text{for }1\le i\le n-1,\ 0\le j\le n-1
\endpmatrix\\
&=
\int\det\pmatrix 
1&u&u^2&\dots&u^{n-1}\\
\mu_0&\mu_1&\mu_2&\dots&\mu_{n-1}\\
\mu_1&\mu_2&\mu_3&\dots&\mu_{n}\\
\hdotsfor5\\
\mu_{n-2}&\mu_{n-1}&\mu_{n+1}&\dots&\mu_{2n-3}
\endpmatrix \frac{d\mu(u)}{u-y_1}\\
&=(-1)^{n-1}\det\limits_{0\le i,j\le n-2}\left(
\mu_{i+j}\right)
\int \frac{p_{n-1}(u)}{u-y_1}\,d\mu(u)\\
&=(-1)^{n}\det\limits_{0\le i,j\le n-2}\left(
\mu_{i+j}\right)
q_{n-1}(y_1),
\endalign
$$
where the next-to-last line is due to Lemma~\TB.
This confirms~(\AB) for $m=0$ and $k=1$.

\medskip
We now turn to the induction step. 
We have to distinguish between
two cases, depending on whether $n\ge k$ or not.
For convenience, we rewrite (\AB) (with~$m=0$) in the form
$$
(-1)^{nk}
\frac {\prod\limits _{1\le i<j\le k} ^{}(y_i-y_j)}
{\det\limits_{0\le i,j\le n-k-1}\left(\mu_{i+j}\right)}
\det\limits_{0\le i,j\le n-1}\left(
\dsize\int u^{i+j}\dfrac {d\mu(u)}
{\prod _{\ell=1} ^{k}(u-y_\ell)}\right)
=
{\displaystyle\det M_{k,0,n}(y_1,\dots,y_k)},
\tag\CB
$$
and we apply a similar rewriting to (\AC),
leading to (\CB) with denominator omitted
on the left-hand side and $M$ replaced by $N$ on the right-hand side.

Let $k\ge2$ and assume that
(\CB), and also the analogue corresponding to~(\AC),
is true for ``smaller~$k$".

We are going to use the condensation formula of Proposition~\TD.
For $n\ge k$, the identity (\BC) with $N=k$, $A=
M_{k,0,n}(y_1,\dots,y_k)$, $i_1=j_1=1$ and $i_k=j_k=k$ gives
$$\multline
\det M_{k,0,n}(y_1,\dots,y_k)\cdot 
\det M_{k-2,0,n-1}(y_2,\dots,y_{k-1})\\
=\det M_{k-1,0,n}(y_2,\dots,y_k)\cdot 
\det M_{k-1,0,n-1}(y_1,\dots,y_{k-1})\\
-
\det M_{k-1,0,n-1}(y_2,\dots,y_k)
\cdot \det M_{k-1,0,n}(y_1,\dots,y_{k-1}).
\endmultline
\tag\CC$$
For $n<k$, with the same choices of $N,i_1,i_2,j_1,j_2$, but with
$A=N_{k,0,n}(y_1,\dots,y_k)$, we obtain the same identity, but
with~$N$ instead of~$M$. A little detail is that, if $n=k-1$,
we encounter the terms $N_{k-2,0,n-1}(\dots)=N_{k-2,0,k-2}(\dots)$ and
$N_{k-1,0,n}(\dots)=N_{k-1,0,k-1}(\dots)$ in~(\CC) (with~$M$ replaced
by~$N$). It can be seen by inspection that
$N_{k-2,0,k-2}(\dots)=M_{k-2,0,k-2}(\dots)$ and
$N_{k-1,0,k-1}(\dots)=M_{k-1,0,k-1}(\dots)$. This is the
interpretation that we give these terms in that special case.

The identity (\CC) can be seen 
as a recurrence formula for $\det M_{k,0,n}(y_1,\dots,y_k)$,
as\linebreak one can use it to express $\det M_{k,0,n}(y_1,\dots,y_k)$ in terms
of expressions of the form\linebreak $\det M_{l,0,s}(y_a,\dots,y_b)$ with~$l$
smaller than~$k$ (and similarly with~$M$ replaced by~$N$),
provided the determinants $\det M_{k-2,0,n-1}(y_2,\dots,y_{k-1})$
and $\det N_{k-2,0,n-1}(y_2,\dots,y_{k-1})$ are all non-zero.
(We shall address the latter point later.) 
Hence, for carrying out the induction step it suffices to
verify that the left-hand side of~(\CB) satisfies the same
recurrence. Consequently, we substitute this left-hand side in~(\CC).
After cancellation of factors that are common to both sides, we
arrive at
$$\multline
(y_k-y_1)
\det\limits_{0\le i,j\le n-1}\left(
\dsize\int u^{i+j}\dfrac {d\mu(u)}
{\prod _{\ell=1} ^{k}(u-y_\ell)}\right)
\cdot 
\det\limits_{0\le i,j\le n-2}\left(
\dsize\int u^{i+j}\dfrac {d\mu(u)}
{\prod _{\ell=2} ^{k-1}(u-y_\ell)}\right)\\
=\det\limits_{0\le i,j\le n-1}\left(
\dsize\int u^{i+j}\dfrac {d\mu(u)}
{\prod _{\ell=2} ^{k}(u-y_\ell)}\right)
\cdot 
\det\limits_{0\le i,j\le n-2}\left(
\dsize\int u^{i+j}\dfrac {d\mu(u)}
{\prod _{\ell=1} ^{k-1}(u-y_\ell)}\right)\\
-
\det\limits_{0\le i,j\le n-2}\left(
\dsize\int u^{i+j}\dfrac {d\mu(u)}
{\prod _{\ell=2} ^{k}(u-y_\ell)}\right)
\cdot \det\limits_{0\le i,j\le n-1}\left(
\dsize\int u^{i+j}\dfrac {d\mu(u)}
{\prod _{\ell=1} ^{k-1}(u-y_\ell)}\right).
\endmultline
$$
This is the special case of Lemma~{\TE} where 
$$c_n=\int u^{i+j}\dfrac {d\mu(u)}
{\prod _{\ell=1} ^{k}(u-y_\ell)}
,$$
$\al=-y_k$ and $\be=-y_1$. Since Lemma~\TGa\ with $m=0$ says that
all these determinants do not vanish identically,
this establishes the induction step and
proves~(\CB) and thus the proposition.\quad \quad \qed
\enddemo

We are now ready to prove Theorem~\TA.

\demo{Proof of Theorem \TA}
We apply induction with respect to $k+m$.
As start of the induction we use Propositions~\TG\ and~\TH.
In other words, we know that Theorem~\TA\ holds for $k=0$ and for
$m=0$.

In preparation of the induction step, we again rewrite (\AB),
$$
\multline
(-1)^{n(m-k)+km}
\frac {\bigg(\prod\limits _{1\le i<j\le m} ^{}(x_j-x_i)\bigg)
\bigg(\prod\limits _{1\le i<j\le k} ^{}(y_i-y_j)\bigg)
}
{\det\limits_{0\le i,j\le n-k-1}\left(\mu_{i+j}\right)}\\
\times
\det\limits_{0\le i,j\le n-1}\left(
\dsize\int u^{i+j}\dfrac {\prod _{\ell=1} ^{m}(u-x_\ell)}
{\prod _{\ell=1} ^{k}(u-y_\ell)}\,d\mu(u)\right)\\
=
{\displaystyle\det M_{k,m,n}(x_1,\dots,x_m,y_1,\dots,y_k)},
\endmultline
\tag\CD
$$
and similarly (\AC), leading to (\CD) with denominator omitted
on the left-hand side and $M$ replaced by $N$ on the right-hand side.

Let now $k$ and $m$ be positive integers, and assume that (\CD),
and also its analogue corresponding to~(\AC), hold for ``smaller $k+m$".

Here again,
we are going to use the condensation formula of Proposition~\TD.
For $n\ge k$, the identity (\BC) with $N=k$, $A=
M_{k,m,n}(x_1,\dots,x_m,y_1,\dots,y_k)$, $i_1=j_1=1$ and $i_k=j_k=k+m$ gives
$$\multline
\det M_{k,m,n}(x_1,\dots,x_m,y_1,\dots,y_k)\cdot 
\det M_{k-1,m-1,n}(x_2,\dots,x_m,y_1,\dots,y_{k-1})\\
=\det M_{k,m-1,n+1}(x_2,\dots,x_m,y_1,\dots,y_k)\cdot 
\det M_{k-1,m,n-1}(x_1,\dots,x_m,y_1,\dots,y_{k-1})\\
-
\det M_{k,m-1,n}(x_2,\dots,x_m,y_1,\dots,y_k)
\cdot \det M_{k-1,m,n}(x_1,\dots,x_m,y_1,\dots,y_{k-1}).
\endmultline
\tag\CE$$
For $n<k$, with the same choices of $N,i_1,i_2,j_1,j_2$, but with
$A=N_{k,m,n}(y_1,\dots,y_k)$, we obtain the same identity, but
with~$N$ instead of~$M$. Again we have to take notice of the little
detail that, if $n=k-1$, we encounter terms such as
$N_{k,m-1,n+1}(\dots)=N_{k,m-1,k}(\dots)$ and
$N_{k-1,m,n}(\dots)=N_{k-1,m,k-1}(\dots)$ in~(\CE) (with~$M$ replaced
by~$N$). Here also they have to be interpreted as the corresponding
``$M$-terms".

In order to accomplish the induction step, we have to prove that the
left-hand side of~(\CD) satisfies the same relation. 
 Consequently, we substitute this left-hand side in~(\CE).
After cancellation of factors that are common to both sides, we
arrive at
$$\multline
\det\limits_{0\le i,j\le n-1}\left(
\dsize\int u^{i+j}\dfrac {\prod _{\ell=1} ^{m}(u-x_\ell)}
{\prod _{\ell=1} ^{k}(u-y_\ell)}\,d\mu(u)\right)\cdot
\det\limits_{0\le i,j\le n-1}\left(
\dsize\int u^{i+j}\dfrac {\prod _{\ell=2} ^{m}(u-x_\ell)}
{\prod _{\ell=1} ^{k-1}(u-y_\ell)}\,d\mu(u)\right)\\
=
-\det\limits_{0\le i,j\le n}\left(
\dsize\int u^{i+j}\dfrac {\prod _{\ell=2} ^{m}(u-x_\ell)}
{\prod _{\ell=1} ^{k}(u-y_\ell)}\,d\mu(u)\right)\cdot
\det\limits_{0\le i,j\le n-2}\left(
\dsize\int u^{i+j}\dfrac {\prod _{\ell=1} ^{m}(u-x_\ell)}
{\prod _{\ell=1} ^{k-1}(u-y_\ell)}\,d\mu(u)\right)
\kern.5cm\\
+
\det\limits_{0\le i,j\le n-1}\left(
\dsize\int u^{i+j}\dfrac {\prod _{\ell=2} ^{m}(u-x_\ell)}
{\prod _{\ell=1} ^{k}(u-y_\ell)}\,d\mu(u)\right)
\cdot 
\det\limits_{0\le i,j\le n-1}\left(
\dsize\int u^{i+j}\dfrac {\prod _{\ell=1} ^{m}(u-x_\ell)}
{\prod _{\ell=1} ^{k-1}(u-y_\ell)}\,d\mu(u)\right).
\endmultline
$$
This is the special case of Lemma~{\TF} where 
$$c_n=\int u^{i+j}\dfrac {\prod _{\ell=2} ^{m}(u-x_\ell)}
{\prod _{\ell=1} ^{k}(u-y_\ell)}\,d\mu(u)
,$$
$\al=-x_1$ and $\be=-y_k$. Since Lemma~\TGa\ says that
all these determinants do not vanish identically,
this establishes the induction step and
proves~(\CD), and thus Theorem~\TA.\quad \quad \qed
\enddemo

\subhead 5. The case of equal parameters\endsubhead
Let
$$
R(x_1,x_2,\dots,x_m,y_1,y_2,\dots,y_k)
$$
denote the right-hand side of (\AB) if $n\ge k$, and
the right-hand side of (\AC) if $n< k$.
Since the numerator (regardless whether we are considering~(\AB)
or~(\AC)) is skew-symmetric in the $x_i$'s and 
skew-symmetric in the $y_i$'s, it is
divisible by the Vandermonde products 
$$\bigg(\prod _{1\le i<j\le m} ^{}(x_j-x_i)\bigg)
\bigg(\prod _{1\le i<j\le k} ^{}(y_i-y_j)\bigg)$$ 
in the denominator. Thus, while in its
definition it seems problematic to substitute the same value
for two different $x_i$'s or for two different $y_i$'s 
in $R(x_1,x_2,\dots,x_m,y_1,y_2,\dots,y_k)$, this is
actually not the case. The proposition below provides an 
explicit expression for such substitutions of equal 
parameters.

\proclaim{Proposition \TI}
Let $r,s$, $k_1,k_2,\dots,k_s$, $m_1,m_2,\dots,m_r$, and $n$ 
be non-negative integers
and $\xi_1,\xi_2,\dots,\xi_r$ and $\om_1,\om_2,\dots,\om_s$ be variables.
We write $k$ for the sum $\sum_{i=1}^sk_i$, and $m$ for the sum
$\sum_{i=1}^rm_i$.

Then
$$
\multline
R(\xi_1,\dots,\xi_1,\xi_2,\dots,\xi_2,\dots,\xi_r,\dots,\xi_r,
\om_1,\dots,\om_1,\om_2,\dots,\om_2,\dots,\om_s,\dots,\om_s)
\\
=(-1)^{n(m-k)+km}
\frac {\displaystyle\det M_{k_1,\dots,k_s,m_1,\dots,m_r,n}
(\xi_1,\dots,\xi_r,\om_1,\dots,\om_s)}
{
\bigg(\prod\limits _{1\le i<j\le r} ^{}(x_j-x_i)^{m_im_j}\bigg)
\bigg(\prod\limits _{1\le i<j\le s} ^{}(y_i-y_j)^{k_ik_j}\bigg)
},
\endmultline
$$
where $\xi_i$ is repeated $m_i$ times 
and $\om_i$ is repeated $k_i$ times 
in the argument of $R$ on the
left-hand side. The matrix in the numerator on the right-hand side
is defined by
$$
M_{k_1,\dots,k_s,m_1,\dots,m_r,n}
(\xi_1,\dots,\xi_r,\om_1,\dots,\om_s)
=
\pmatrix P_{m_1,k,n}(\xi_1)\\
\vdots\\
P_{m_r,k,n}(\xi_r)\\
Q_{k_1,k,n}(\om_1)\\
\vdots\\
Q_{k_s,k,n}(\om_s)
\endpmatrix,
$$
with
$$
P_{a,k.n}(\xi)=\pmatrix 
\dfrac {p^{(i-1)}_{n-k+j-1}(\xi)} {(i-1)!}
\endpmatrix_{1\le i\le a,\ 1\le j\le k+m}
$$
and
$$
Q_{a,k,n}(\om)=\pmatrix 
\dfrac {q^{(i-1)}_{n-k+j-1}(\om)} {(i-1)!}
\endpmatrix_{1\le i\le a,\ 1\le j\le k+m}.
$$
If $b<0$ the polynomial $p_b(\xi)$ has to be interpreted as~$0$,
while the function $q_b(\om)$ has to be interpreted as~$\om^{-b-1}$,
\endproclaim

\demo{Proof}
This can be proved in the same way as \cite{\KratCN, Prop.~5}.
We leave the details to the reader.\quad \quad \qed
\enddemo

\remark{Remark}
Obviously, we have
$$
q^{(i-1)}_n(y)=(-1)^{i-1}(i-1)!\int \frac {p_n^{(i-1)}(u)} {(y-u)^i}\,d\mu(u).
$$
\endremark

\subhead 6. Uvarov's formula\endsubhead
Uvarov's theorem \cite{\UvarAA, \UvarAB} (cf\. \cite{\IsmaAA,
Theorem~2.7.3}) says that, {\it in the setting of Theorem~{\rm\TA},
the right-hand sides of\/}~(\AB) {\it and\/}~(\AC),
{\it seen as polynomials in $x_1$, give orthogonal
polynomials for the density}
$$
\dfrac {\prod _{\ell=2} ^{m}(u-x_\ell)}
{\prod _{\ell=1} ^{k}(u-y_\ell)}\,d\mu(u).
\tag\DA
$$
(The reader should note that the product in the numerator starts with
$\ell=2$.)

The connection with our Theorem~\TA\ is set up by Lemma~\TC.
Namely, if in Lemma~\TC\ we choose for the moments $\nu_i$ the moments
corresponding to the density in (\DA), then the determinant
(\BBf) turns out to be exactly the determinant on the left-hand sides
of~(\AB) and~(\AC). Hence, in view of Lemma~\TC, it is obvious that
the right-hand sides of~(\AB) and~(\AC) are orthogonal polynomials
with respect to the density in~(\DA). There is one ``catch" however:
Lemma~\TC\ says that this is only the case if all Hankel determinants
of moments of~(\DA) are non-zero. Phrased differently: the determinant
in (\BBf) must be a polynomial in~$x$ of degree~$n$, for $n=0,1,\dots$.

Uvarov's theorem provides one such scenario: he shows that, if
$$
\dfrac {\prod _{\ell=2} ^{m}(u-x_\ell)}
{\prod _{\ell=1} ^{k}(u-y_\ell)}
$$
is positive for all $u$ in the support of the measure $\mu(u)$,
then the right-hand sides of~(\AB) and~(\AC) are polynomials in~$x_1$
of degree~$n$.

On the other hand, more generally, our Theorem~\TA\ implies that,
{\it whenever the right-hand sides of\/} (\AB) {\it and\/} (\AC) {\it are
polynomials in~$x_1$ of degree~$n$, then they are orthogonal
polynomials in~$x_1$ for the density in} (\DA).

\subhead 7. Applications: some Hankel determinant evaluations\endsubhead
In this final section, we apply Theorem~\TA\ to derive some Hankel
determinant evaluations featuring Catalan numbers and central
binomial coefficients.

\medskip
We restrict ourselves to the special case of Theorem~\TA\ where
$$
d\mu(u)=\sqrt{1-\frac {u^2}4}\,du,  \quad \quad -2\le u\le 2.
\tag\ZA
$$
The polynomials which are orthogonal with
respect to this density are, essentially, 
the {\it Chebyshev polynomials of the
second kind\/} $U_n(x)$, given by the generating function
$$
\sum_{n\ge0}U_n(x)\,z^n=\frac {1} {1-2xz+z^2}.
$$
More precisely, we have
$$
\frac {1} {\pi}\int_{-2}^2U_m(u/2)U_n(u/2)\sqrt{1-\tfrac {u^2}4}\,du=\de_{n,m},
$$
with, again, $\de_{n,m}$ denoting the Kronecker delta.
The polynomials $U_n(x/2)$ are monic, as can be seen from
$$
\sum_{n=0}^\infty U_n(x/2)z^n=\frac {1} {1-xz+z^2}.
\tag\ZB$$
The non-vanishing 
moments of the density (\ZA) are given by the {\it Catalan numbers}
$C_n=\frac {1} {n+1}\binom {2n}n$. To be precise, we have
$$
\frac {1} {\pi}\int_{-2}^2 u^n\sqrt{1-\frac {u^2}4}\,du
=\cases C_{n/2},&\text{if $n$ is even},\\
0,&\text{if $n$ is odd},
\endcases\quad \quad n\ge0.
$$

\medskip
We shall apply Theorem~\TA\ with $k=1$ and $m=0$, and with
$k=m=1$. For carrying out the corresponding calculations, we
need the following fundamental integral evaluation:\footnote{The 
substitution $u=2t/(1+t^2)$ transforms the
integral into an integral of a rational function. Subsequently, a
thorough case analysis has to be performed.}
$$
\frac {1} {\pi}\int_{-1}^1\frac {\sqrt{1-u^2}} {u+a}\,du=
\cases 
a-\sqrt{a^2-1},&\text{for }a\ge1,\\
&\text{for }\Re a>0\text{ and }\Im a\ne0,\\
&\text{for }\Re a=0\text{ and }\Im a>0,\\
a,&\text{for }-1\le a\le 1,\\
a+\sqrt{a^2-1},&\text{for }a\le-1,\\
&\text{for }\Re a<0\text{ and }\Im a\ne0,\\
&\text{for }\Re a=0\text{ and }\Im a<0.
\endcases
\tag\ZC
$$

For a ``test", we may use this to confirm the orthogonality
of the $U_n(x/2)$'s. We have
$$\align
\sum_{n\ge0}z^n\frac {1} {\pi}\int_{-2}^2 U_n(u/2)\sqrt{1-\tfrac {u^2}4}\,du
&=
\frac {1} {\pi}\int_{-2}^2\frac {\sqrt{1-\tfrac {u^2}4}}
{1-uz+z^2}\,du\\
&=
\frac {2} {\pi}\int_{-1}^1\frac {\sqrt{1-u^2}} {1-2uz+z^2}\,du\\
&=
-\frac {1} {z\pi}\int_{-1}^1\frac {\sqrt{1-u^2}} {u-\frac {1+z^2} {2z}}\,du.
\tag\ZD
\endalign$$
This calculation is meant as a calculation for (formal) power series
in~$z$. Thus, we must think of~$z$ as ``small", which implies that
$\frac {1+z^2} {2z}$ is ``large". Hence, when we apply (\ZC) with $a=-\frac
{1+z^2} {2z}$, we find ourselves in the third case of the case
distinction. This gives
$$
\sum_{n\ge0}z^n\frac {1} {\pi}\int_{-2}^2 U_n(u/2)\sqrt{1-\tfrac
{u^2}4}\,du
=-\frac {1} {z}\left(-\frac {1+z^2}{2z}+\sqrt{\frac {(1+z^2)^2} {4z^2}-1}\right)
=1,
\tag\ZE
$$
as expected.
Another ``test" that we may perform is the calculation of moments
(of the measure) of the Chebyshev polynomials.
We have
$$\align
\sum_{n\ge0}z^n\frac {1} {\pi}\int_{-2}^2u^n\sqrt{1-\tfrac {u^2}4}\,du&=
\frac {1} {\pi}\int_{-2}^2\frac {\sqrt{1-\tfrac {u^2}4}} {1-uz}\,du
=
\frac {2} {\pi}\int_{-1}^1\frac {\sqrt{1-u^2}} {1-2uz}\,du\\
&=
-\frac {1} {z\pi}\int_{-1}^1\frac {\sqrt{1-u^2}} {u-\frac {1} {2z}}\,du
=
\frac {1} {2z^2}\left(1-\sqrt{1-4z^2}\right)\\
&=\sum_{n\ge0}C_nz^{2n},
\tag\ZF
\endalign$$
again as expected.

We intend to apply Theorem~\TA\ with $k=1$ and $y_1=-2a$.
For convenience, let us write $\om(a)$ for the quantity in~(\ZC).
Then, by~(\ZC) and~(\ZF), 
for the moments of the density $\frac {1} {u-y_1}\sqrt{1-\frac
{u^2} {4}}\,du$ we have
$$\align 
\sum_{n\ge0}z^n\frac {1} {\pi}&
\int_{-2}^2\frac {u^n\sqrt{1-\frac {u^2}4}} {u+2a} \,du
=
\frac {1} {\pi}
\int_{-2}^2\frac {\sqrt{1-\frac {u^2}4}} {(u+2a)(1-uz)} \,du\\
&=
\frac {1} {\pi(1+2az)}
\int_{-2}^2\frac {\sqrt{1-\frac {u^2}4}} {u+2a} \,du
+
\frac {z} {\pi(1+2az)}
\int_{-2}^2\frac {\sqrt{1-\frac {u^2}4}} {1-uz} \,du\\
&=
\frac {\om(a)} {1+2az}
+
\frac {z} {1+2az}\sum_{n\ge0}C_nz^{2n}\\
&=
\sum_{n\ge0}\om(a)(-2a)^nz^n
+\sum_{n\ge0}z^n
\sum_{k=0}^{\fl{(n-1)/2}}(-2a)^{n-2k-1}
C_k.
\endalign$$
Hence,
$$
\frac {1} {\pi}
\int_{-2}^2\frac {u^n\sqrt{1-\frac {u^2}4}} {u+2a} \,du
=
\om(a)(-2a)^n+
\sum_{k=0}^{\fl{(n-1)/2}}(-2a)^{n-2k-1}
C_k.
\tag\ZG
$$
We next compute
the functions $q_n(-2a)$ associated with the Chebyshev polynomials\linebreak
$U_n(x/2)$ (cf\. (\AA)),
$$\align 
\sum_{n\ge0}z^nq_n(-2a)
&=\sum_{n\ge0}z^n\frac {1} {\pi}
\int_{-2}^2\frac {U_n(u/2)\sqrt{1-\frac {u^2}4}} {-2a-u} \,du
=
-\frac {1} {\pi}
\int_{-2}^2\frac {\sqrt{1-\frac {u^2}4}} {(u+2a)(1-uz+z^2)} \,du\\
&=
-\frac {1} {\pi(1+2az+z^2)}
\int_{-2}^2\frac {\sqrt{1-\frac {u^2}4}} {u+2a} \,du\\
&\kern3cm
-
\frac {z} {\pi(1+2az+z^2)}
\int_{-2}^2\frac {\sqrt{1-\frac {u^2}4}} {1-uz+z^2} \,du\\
&=
-\frac {\om(a)+z} {1+2az+z^2}
=-\sum_{n\ge0}\big(\om(a)U_n(-a)+U_{n-1}(-a)\big)z^n,
\endalign$$
where we used (\ZC), (\ZD)/(\ZE), and finally~(\ZB) to obtain the last line.
Hence,
$$
q_n(-2a)=-\big(\om(a)U_n(-a)+U_{n-1}(-a)\big).
\tag\ZH$$

If we put these findings together, we get the following result
from Theorem~\TA\ with $k=1$ and $m=0$.

\proclaim{Theorem \TJ}
Let $X$ and $a$ be variables.
For all positive integers $n$, we have
$$\multline
\det_{0\le i,j\le n-1}\left(X(-2a)^{i+j}
+\sum_{k=0}^{\fl{(i+j-1)/2}}(-2a)^{i+j-2k-1}C_k\right)\\
=(-1)^{n-1}\big(XU_{n-1}(-a)+U_{n-2}(-a)\big).
\endmultline\tag\ZI$$
\endproclaim

\demo{Proof}
Using (\ZG) and (\ZH), we see that
Theorem~\TA\ with $k=1$, $m=0$, $y_1=-2a$, and $d\mu(u)=\sqrt{1-\frac
{u^2} {4}}\,du$ reads
$$\multline
\det_{0\le i,j\le n-1}\left(\om(a)(-2a)^{i+j}
+\sum_{k=0}^{\fl{(i+j-1)/2}}(-2a)^{i+j-2k-1}C_k\right)\\
=(-1)^{n-1}\big(\om(a)U_{n-1}(-a)+U_{n-2}(-a)\big).
\endmultline$$
This holds for all complex numbers $a$. Now, 
if we regard the above determinant (formally) as a polynomial
in $\om(a)$, then it is not difficult to see that the degree
in~$\om(a)$ is at most~1. On the other hand, the right-hand side
of the above equation is visibly a polynomial of degree~1 
in~$\om(a)$. If we now choose $a=-2$ (say), so that
$a+\sqrt{a^2-1}=-2+\sqrt3$ and $1$ are linearly independent over the
rational numbers, then we may conclude that $\om(a)$ can be replaced by
a variable, $X$ say. This establishes the assertion of the
theorem.\quad \quad \qed
\enddemo

If we choose $a=-1$ in the above theorem, the sum in the determinant
can be simplified. Namely, we have
$$
(-1)2^n
+\sum_{k=0}^{\fl{(n-1)/2}}2^{n-2k-1}C_k
=\cases -\binom {n}{n/2},&\text{if $n$ is even},\\
-\frac {1} {2}\binom {n+1}{(n+1)/2},&\text{if $n$ is odd},
\endcases
\tag\ZIa$$
as is straightforward to verify by an induction on~$n$.
Since $U_n(1)=n+1$ for  all~$n$, the choice of~$X=-1$ in Theorem~\TJ\ leads
to the Hankel determinant evaluation
$$
\det_{0\le i,j\le n-1}\left(2^{-2\cl{(i+j)/2}}
\binom {2\cl{(i+j)/2}} {\cl{(i+j)/2}}\right)
=2^{-n(n-1)}.
\tag\ZJ$$
More generally, if we choose $X=-Y-1$, then we get
$$
\det_{0\le i,j\le n-1}\left(Y+
2^{-2\cl{(i+j)/2}}
\binom {2\cl{(i+j)/2}} {\cl{(i+j)/2}}\right)
=2^{-n(n-1)}(Yn+1).
\tag\ZK$$

In a similar fashion, Theorem~\TA\ with $k=m=1$, $x_1=b$, $y_1=-2a$,
and $d\mu(u)=\sqrt{1-\frac {u^2} {4}}\,du$ implies the following result.

\proclaim{Theorem \TK}
Let $X$, $a$ and $b$ be variables.
For all positive integers $n$, we have
$$\multline
\det_{0\le i,j\le n-1}\left(
X(-2a)^{i+j+1}
+\sum_{k=0}^{\fl{(i+j)/2}}(-2a)^{i+j-2k}C_k\right.\\
\left.
-b\bigg(X(-2a)^{i+j}
+\sum_{k=0}^{\fl{(i+j-1)/2}}(-2a)^{i+j-2k-1}C_k\bigg)
\right)\\
=U_{n-1}(b/2)\big(XU_{n}(-a)+U_{n-1}(-a)\big)
-U_{n}(b/2)\big(XU_{n-1}(-a)+U_{n-2}(-a)\big).
\endmultline\tag\ZL$$
\endproclaim

Again, it is worth stating the identities that one obtains
in the special case where $a=-1$ explicitly.
Namely, in view of the simplification~(\ZIa) and of $U_n(1)=n+1$, 
in that case Equation~(\ZL) reduces to
%
%
$$\multline
\det_{0\le i,j\le n-1}\left(
Y
+2^{-2\cl{(i+j+1)/2}}\binom {2\cl{(i+j+1)/2}} {\cl{(i+j+1)/2}}
\right.\\
\left.
-b\bigg(\tfrac {1} {2}Y
+2^{-2\cl{(i+j)/2}-1}\binom {2\cl{(i+j)/2}} {\cl{(i+j)/2}}
\bigg)
\right)\\
=(-1)^{n-1}2^{-n^2}\Big(U_{n-1}(b/2)(Y(n+1)+1)
-U_{n}(b/2)(Yn+1)\Big).
\endmultline\tag\ZM$$
Since $U_n(0)=(-1)^{n/2}$ if $n$ is even and $U_n(0)=0$ otherwise,
for $b=0$ Equation~(\ZM) becomes
$$
\det_{0\le i,j\le n-1}\left(
Y
+2^{-2\cl{(i+j+1)/2}}\binom {2\cl{(i+j+1)/2}} {\cl{(i+j+1)/2}}
\bigg)
\right)
=(-1)^{\binom n2}2^{-n^2}\big(2\cl{n/2}Y+1\big).
\tag\ZN$$
Finally, since
$$
U_n(1/2)=\cases 
1,&\text{if }n\equiv 0,1~\text{mod }6,\\
0,&\text{if }n\equiv 2~\text{mod }3,\\
-1,&\text{if }n\equiv 3,4~\text{mod }6,
\endcases
$$
for $b=1$ Equation (\ZM) becomes
$$\multline
\det_{0\le i,j\le n-1}\left(
Y
+2^{-2\cl{(i+j)/2}}\binom {2\cl{(i+j)/2}} {\cl{(i+j+1)/2}}
\bigg)
\right)\\
=\cases 
2^{-n(n-1)}Y,&\text{if }n\equiv0~(\text{mod 3}),\\
-2^{-n(n-1)}(Y(n+1)+1),&\text{if }n\equiv1~(\text{mod 3}),\\
2^{-n(n-1)}(Yn+1),&\text{if }n\equiv2~(\text{mod 3}).
\endcases
\endmultline\tag\ZO$$

When testing (\ZK) and (\ZN) with the help of
{\sl Mathematica}, initially I mistyped the
exponent of~2. This led to the discovery of
two more --- conjectural --- Hankel determinant evaluations.

\proclaim{Conjecture \TL}
Let $Y$ be a variable.
For all positive integers $n$, we have
$$
\det_{0\le i,j\le n-1}\left(Y+
2^{-i-j}
\binom {2\cl{(i+j)/2}} {\cl{(i+j)/2}}\right)
=-2^{-(n-1)^2}\big(Y(n-3)+1\big)
\tag\ZP
$$
and
$$\multline
\det_{0\le i,j\le n-1}\left(
Y
+2^{-i-j-1}\binom {2\cl{(i+j+1)/2}} {\cl{(i+j+1)/2}}
\bigg)
\right)\\
=\cases
(-1)^{\fl{n/6}}2^{-n(n-1)}
\left(\frac {1} {2}\left\lfloor\frac {4n+2} {3}\right\rfloor
Y+1\right),&\text{if $n$ is even},\\
(-1)^{\fl{n/6}}2^{-n(n-1)}
\left(\frac {1} {2}\left\lfloor\frac {4n+4} {3}\right\rfloor
Y+1\right),&\text{if $n$ is odd}.
\endcases
\endmultline\tag\ZQ$$
\endproclaim

\remark{Acknowledgement} 
I thank Mourad Ismail for insisting that, after having
written~\cite{\KratCN}, I throw a ``combinatorial
eye" on  Uvarov's formula.
\endremark

\enddocument